\documentclass[12pt]{amsart}
\usepackage{amsfonts}
\usepackage{enumerate}
\usepackage{helvet,courier}
\usepackage{mathptmx,amsmath}
\usepackage{type1cm}
\DeclareMathAlphabet{\mathcal}{OMS}{cmsy}{m}{n}
\DeclareMathAlphabet{\mathbbold}{U}{bbold}{m}{n}  
\usepackage{amsthm}
\usepackage{graphicx}       
\usepackage{multicol}
\usepackage{mathrsfs}
\usepackage{amssymb}
\usepackage{verbatim}

\usepackage[usenames,dvipsnames]{xcolor}

\theoremstyle{plain}
\newtheorem{thm}{Theorem}[section]
\newtheorem{lm}[thm]{Lemma}

\newtheorem{prop}[thm]{Proposition}

\theoremstyle{remark}
\newtheorem{rmk}{Remark}

\theoremstyle{definition}

\newcommand{\bnu}{\begin{enumerate}}
\newcommand{\enu}{\end{enumerate}}
\newcommand{\q}{\quad}
\newcommand{\qq}{\qquad}

\newcommand{\sset}{\subset}

\newcommand{\ga}{\gamma}

\newcommand{\om}{\omega}
\newcommand{\Om}{\Omega}

\newcommand{\la}{\lambda}

\newcommand{\si}{\sigma}

\newcommand{\de}{\delta}

\newcommand{\bbz}{\mathbb{Z}}

\newcommand{\bbzn}{\mathbb Z^n}

\newcommand{\bbr}{\mathbb{R}}

\newcommand{\bbrn}{\mathbb R^n}
\newcommand{\rn}{\mathbb{R}^n}

\newcommand{\bbn}{\mathbb{N}}

\newcommand{\UU}{\mathcal{U}}
\newcommand{\VV}{\mathcal{V}}
\newcommand{\II}{\mathcal{I}}
\newcommand{\TT}{\mathcal{T}}
\newcommand{\LL}{\mathcal{L}}

\newcommand{\kkk}{\vec{\boldsymbol{k}}}
\newcommand{\xxxi}{\vec{\boldsymbol{\xi}}}
\newcommand{\xxx}{\vec{\boldsymbol{x}}}
\newcommand{\yyy}{\vec{\boldsymbol{y}}}

\newcommand{\GGG}{\vec{\boldsymbol{G}}}
\newcommand{\0}{\vec{\boldsymbol{0}}}

\newcommand{\f}{\frac}

\newcommand{\nf}{\infty}

\newcommand{\tf}{\tfrac}
\newcommand{\wh}{\widehat}

\newcounter{question}

\newcommand{\bpf}{\begin{proof}}

\newcommand{\epf}{\end{proof}}

\allowdisplaybreaks

\makeindex         

\begin{document}

\author{Loukas Grafakos}
\address{Department of Mathematics, University of Missouri, Columbia, MO 65211, USA} 
\email{grafakosl@missouri.edu}

\author{Danqing He}
\address{School of Mathematical Sciences,
Fudan University, People's Republic of China}
\email{hedanqing@fudan.edu.cn}

\author{Petr Honz\'ik}
\address{Department of Mathematics,
Charles University, 116 36 Praha 1, Czech Republic}
\email{honzik@gmail.com}

\author{Bae Jun Park}
\address{School of Mathematics, Korea Institute for Advanced Study, Seoul 02455, Republic of Korea}
\email{qkrqowns@kias.re.kr}

\thanks{L. Grafakos would like to acknowledge the support of  the Simons Foundation grant 624733. 
D. He is supported by   NNSF of China (No. 11701583).
P. Honz\'ik was supported by  GA\v CR P201/18-07996S.
B. Park is supported in part by NRF grant 2019R1F1A1044075 and by a KIAS Individual Grant MG070001 at the Korea Institute for Advanced Study} 

\title{Initial $L^2\times\cdots\times L^2 $ bounds for multilinear operators}

\begin{abstract} 
The $L^p$ boundedness theory of convolution operators is  \linebreak based on
 an initial $L^2\to L^2$ estimate    derived from the Fourier transform.   
The corresponding theory of multilinear operators lacks 
such a simple initial estimate in view  of the unavailability of 
Plancherel's identity in this setting, and up to now it
has not been clear  what a natural initial estimate might be. 
In this work we 
%achieve exactly this goal, i.e.,   
obtain  initial    $L^2\times\cdots\times L^2\to L^{2/m}$ estimates   for three types of important multilinear operators: 
 rough  singular integrals, multipliers of H\"ormander type, and multipliers whose derivatives satisfy qualitative estimates. 
%Our main tools are compactly supported wavelets and a combinatorial argument applied to   general 
%building blocks of $m$-linear multiplier operators. 

\end{abstract}

\maketitle

%\tableofcontents

\section{Introduction and Preliminaries}

The systematic study of multilinear operators in  harmonic analysis was initiated by Coifman and Meyer in the seventies. Many important  multilinear operators  can be   written as
\[
T(f_1,\dots , f_m)(x) =   W* (f_1\otimes \cdots \otimes f_m)(x,\dots, x), \qquad x\in \bbrn,
\]
where $f_j$ are Schwartz functions on $\bbrn$,  $W$ is a tempered distribution on $(\bbrn)^m$, and 
$(f_1\otimes \cdots \otimes f_m)(x_1,\dots , x_m) = 
f_1(x_1)\cdots f_m(x_m)$.  Alternatively $T(f_1,\dots , f_m)(x)$ can be expressed as 
\begin{equation}\label{FTofW}
\int_{\bbrn} \cdots \int_{\bbrn} \widehat{f_1}(\xi_1)\cdots \widehat{f_m}(\xi_m) \widehat{W}(\xi_1,\dots , \xi_m)
e^{2\pi i \langle x  , \xi_1+\cdots +\xi_m\rangle } d\xi_1\cdots d\xi_n, 
\end{equation}
where $\widehat{f_j}(\xi)= \int_{\bbrn} f_j(x)e^{-2\pi i\langle x,\xi\rangle}dx$
  denotes the Fourier transform of  $f_j$ and $\widehat{W}$ is the distributional 
Fourier transform of $W$, which must be an $L^\nf$ function if $T$ is bounded 
from $L^{p_1}(\rn)\times \cdots \times L^{p_m}(\rn)$ to $L^{p }(\rn)$ for some choice of indices 
that satisfy $1/p=1/p_1+\cdots+1/p_m$. 
If such an estimate holds for $T$, then $\widehat{W}$ is called a multilinear Fourier multiplier.    
The first important result concerning multilinear Fourier multipliers is  
a nontrivial adaptation of Mihlin's multiplier  condition,  obtained by  Coifman and Meyer    \cite{CM}. The  proof they gave is based on   decomposing the multiplier as a sum of products 
of linear operators via Littlewood-Paley and 
Fourier series expansions. This powerful idea       has essentially been the only technique
available  in this area until the appearance of the wave-packet decompositions
in the work of   Lacey and Thiele \cite{LT1,LT2} on the bilinear Hilbert transform.

If the distribution $W$ has the form 
\[
W=\textup{p.v.}  \,\, \f{1 }{|(y_1,\dots , y_m)|^{ mn}} \Omega \Big( \f{(y_1,\dots , y_m)}{|(y_1,\dots , y_m)|}  \Big)
\]
for some integrable function $\Omega$ on the sphere $\mathbb S^{mn-1}$ with integral zero, then $T$ is called 
an $m$-linear homogeneous singular. The associated operator 
is bounded  if $\Om$ is smooth but it could be 
unbounded if $\Om$ is merely integrable; see \cite{GHS-CRAS}. In this paper we focus on the intermediate 
situation where $\Om$ lies in $L^q$  for some $q\in (1 , \nf]$; these $\Om$'s give rise to rough 
$m$-linear  homogeneous singular integrals. 
The study of bilinear homogeneous singular integrals was initiated by Coifman and Meyer  
in \cite{CM2} who   addressed the case where $\Om$ is a function of bounded variation. 
The boundedness of $T$ in the 
more difficult case when $\Om $ is merely in $L^\nf$ was not solved until four decades later    in  \cite{Gr_He_Ho} in terms
of wavelet decompositions. Prior to that, the first author and Torres \cite{GT-AIM} had 
proved boundedness for $T$ for any  $m$  when  $\Om$   satisfies  
a Lipschitz condition on the sphere. 
In the case $m=1$ the  known results are much better. For instance, Calder\'on and Zygmund \cite{Ca_Zy}   showed that $T$ is bounded in $L^p(\bbrn)$ for $1<p<\infty$ if $\Omega\in L\log L(\mathbb{S}^{n-1})$. This result was improved by Coifman and Weiss \cite{Co_We} under the less restrictive condition that $\Omega$ belongs to the Hardy space $H^1(\mathbb{S}^{n-1})$.

One fundamental difference between linear convolution operators and   multilinear  convolution operators of   type \eqref{FTofW} is an initial estimate. In the linear case the initial estimate is usually $L^2\to L^2$
and this is obtained by Plancherel's identity. There is   not a 
straightforward initial estimate 
for multilinear operators and in most times, it is difficult to find one. 
Inspired by \cite{Gr_He_Ho}, the first two authors and Slav\'ikov\'a \cite{Gr_He_Sl} obtained a bilinear substitute of the Plancherel criterion for $L^2\times L^2\to L^1$ boundedness    for  multipliers in $L^q(\rn)$ ($0<q<4$) with sufficiently many   bounded derivatives. 
This result 
has also been  proved by   Kato,  Miyachi, and Tomita~\cite{Katoetal}  
and  has found many applications; see  for instance \cite{Gr_He_Sl,PS}.

%In the trilinear case a good initial estimate should be 
%$L^2\times L^2\times L^2 \to L^{2/3}$, but the space $L^{2/3}$ is not at all 
%similar to $L^1$, for instance, it is not locally convex. 
%However, the lack of local convexity does not present as many obstacles in our analysis as the multilinear complexity as $m$ increases. 

Overcoming the combinatorial complexity  that arises
from   multilinearity,   in this paper  we  develop  a  method that yields the crucial 
$$
\overbrace{L^2(\mathbb R^n) \times\cdots\times L^2(\mathbb R^n)}^{\text{$m$ times}}\to L^{2/m}(\mathbb R^n)
$$
  estimates %for a general class of building blocks of $m$-linear operators. Our result  allows us to obtain similar boundedness 
for a variety of $m$-linear operators, including rough homogeneous singular integrals 
and multipliers.  Our results    contribute to  the recent surge of activity in the theory of  rough
multilinear singular integrals, see for instance 
\cite{CoCuDiPliOu,DiPliLiVuMa,DiPliLiVuMa2,Gr_He_Sl,6author,Katoetal}.

\medskip

We  first present a sharp $L^2\times \cdots\times L^2\to L^{2/m}$ boundedness criterion for a  multiplier with bounded derivatives up to a certain order. This provides   a multilinear extension of the main result in  \cite{Gr_He_Sl}.

\begin{thm}\label{application1}
Let $m$ be a positive integer with $m\ge 2$ and $1< q<\frac{2m}{m-1}$. Set $M_q$ to be a positive integer satisfying 
$$
M_q>\frac{m(m-1)n}{2m-(m-1)q}.
$$
Suppose that $\si\in L^q((\bbrn)^{m})\bigcap \mathscr{C}^{M_q}((\bbrn)^{m})$ with 
$$\big\Vert \partial^{\alpha}\si \big\Vert_{L^{\infty}((\bbrn)^{m})}\le D_0, \q \text{for }~ |\alpha|\le M_q.$$
Then the estimate 
\begin{equation*}
\big\Vert T_{\si}(f_1,\dots,f_m)\big\Vert_{L^{2/m}(\bbrn)}\lesssim D_0^{1-\frac{(m-1)q}{2m}}
\Big( \Vert \si\Vert_{L^q((\bbrn)^{m})} \Big)^{\frac{(m-1)q}{2m}}\prod_{j=1}^{m}\Vert f_j\Vert_{L^2(\bbrn)}
\end{equation*}
is valid for Schwartz functions $f_1,\dots,f_m$ on $\bbrn$.
\end{thm}

\hfill

Next we discuss multilinear rough singular integral operators.
For a fixed function   $\Om$  on the unit sphere $\mathbb S^{mn-1}$  and for
 $\yyy':=\yyy/|\yyy|\in \mathbb{S}^{mn-1}$ we let  
\begin{equation}\label{defK}
K(\yyy):=\frac{\Omega(\yyy')}{|\yyy|^{mn}} .
\end{equation} 
We then define the corresponding multilinear operator
\begin{equation*}
\LL_{\Om}\big(f_1,\dots,f_m\big)(x):=p.v. \int_{(\bbrn)^m}{K(\yyy)\prod_{j=1}^{m}f_j(x-y_j)}~d\, \yyy, \q x\in \bbrn
\end{equation*}
 for Schwartz functions $f_1,\dots,f_m$ on $\bbrn$.

\begin{thm}\label{application2} 
 Suppose that $\frac{2m}{m+1}<q\le \infty$ and let $\Om \in L^q(\mathbb{S}^{mn-1})$ 
 satisfying $\int_{\mathbb{S}^{mn-1}} \Om\, d\si_{mn-1}=0$. 
Then there exists a constant $C=C_{n,m,q}>0$ such that
\begin{equation*}
\big\Vert \LL_{\Om}(f_1,\dots,f_m)\big\Vert_{L^{2/m}(\bbrn)}\le C\Vert \Om\Vert_{L^q(\mathbb{S}^{mn-1})}\prod_{j=1}^{m}\Vert f_j\Vert_{L^2(\bbrn)}
\end{equation*}
for Schwartz functions $f_1,\dots,f_m$ on $\bbrn$.
\end{thm}

\smallskip

The last result of this paper is about the boundedness of multilinear multiplier operators of H\"ormander type.
%We also study the boundedness of multilinear multiplier operators of H\"ormander  type. 
The multilinear multiplier operator associated with a bounded function $\si$ on $(\bbrn)^m$ is defined as in (\ref{multiplieroperator});
\begin{equation*}
T_{\si}\big(f_1,\dots,f_m\big)(x):=\int_{(\bbrn)^m}{\si(\xxxi)\Big(\prod_{j=1}^{m}\wh{f_j}(\xi_j)\Big)e^{2\pi i\langle x, \sum_{j=1}^{m}\xi_j\rangle}}d\xxxi
\end{equation*}
for Schwartz functions $f_1,\dots,f_m$ on $\bbrn$.
We choose a Schwartz function $\Phi^{(m)}$ on $(\bbrn)^m$ having the properties that $\wh{\Phi^{(m)}}$ is positive and supported in the annulus $\{\xxxi\in (\bbrn)^m: 1/2\le |\xxxi|\le 2\}$, and $\sum_{\ga\in\bbz}{\wh{\Phi^{(m)}}(\xxxi/2^{\ga})}=1$ for $\xxxi\not= \0$.
In the linear case $m=1$, under the assumption 
$$
\sup_{j\in\bbz}\big\Vert \si(2^j\cdot)\wh{\Phi^{(1)}}\big\Vert_{L_s^q(\bbrn)}<\infty,
$$ 
the condition
$$
s>\max\big( |n/2-n/p|, n/q \big)
$$
implies the boundedness of $T_{\si}$ from $L^p(\bbrn)$ to itself. 
Recently, the bilinear analogue of this result was obtained in the series of papers 
\cite{Gr_He_Ho2018, Gr_Ng2019, Mi_Tom} 
by Grafakos, He, Honz\'ik, Miyachi, Nguyen, and Tomita; all of these results were inspired by the fundamental work of Tomita \cite{To} in this direction.

%We pursue the multilinear extension of these results in the local $L^2$ case.
\begin{thm}\label{application3}
Let $1<q<\infty$ and 
\begin{equation}\label{minimalcondition}
s>\max{((m-1)n/2,mn/q)}.
\end{equation} 
Then there exists an absolute constant $C=C_{n,m,q,s}>0$ such that
\begin{equation*}
\|T_\si(f_1,\dots,f_m)\|_{L^{2/m}(\bbrn)}\le C\sup_{j\in\bbz}\big\Vert \si(2^j \; \vec{\cdot}\; )\wh{\Phi^{(m)}}\big\Vert_{L^q_s((\bbrn)^m)} \prod_{j=1}^{m}\Vert f_j\Vert_{L^2(\bbrn)}
\end{equation*}
for Schwartz functions $f_1,\dots,f_m$ on $\bbrn$.
\end{thm}
We remark that for $1<q\le 2$ this result  has been obtained by \cite{To} and \cite{GS}, so Theorem~\ref{application3} is new 
only in the case $q>2$;  this corresponds to the classical result of Calder\'on and Torchinsky \cite{CT}  in the linear setting.
The sharpness of   condition (\ref{minimalcondition}) was 
addressed in \cite[Theorem 2]{Gr_He_Ho2018}.

\medskip

We design two novel ideas to deal with the above results:
%There are two novel ideas in this paper: 
(I)  An innovative     decomposition of an $m$-linear multiplier into sums of products 
so that  $l$-linear Plancherel type estimates    ($1\le l\le m$) can be used;  see  Proposition~\ref{keyapplication1}. (II)  An effective way  to split lattice points in $(\mathbb Z^{n})^m$
{\it   into groups of columns}   for the purposes of obtaining 
 $L^2\times\cdots\times L^2\to L^{2/m}$ estimates; for  details see  Section~\ref{pfmain}.

It is inevitable to introduce complicated notation  in order to 
comprehensively  present the proofs in the general framework of $m$-linear operators; for this
reason 
we   urge the readers to restrict attention to the case $m=3$, which 
already  presents 
several new ingredients compared with
 the case $m=2$ and contains the main ideas.

\subsection*{Notation}
 $C$ will denote inessential constants that may vary from occurrence to occurrence. $A \lesssim B$
means $A\le CB$ with $C$ independent of $A$ and $B$, and  write $A\approx B$ if both $A \lesssim B$ 
and $B \lesssim A$ hold.
We denote the set of natural numbers by 
 $\bbn$ and of integers by $\bbz$; we also set $\bbn_0:=\bbn\cup \{0\}$. Throughout this paper, 
the index  $m\in  \bbn$ will be   the degree of multilinearity of operators.

\section{Plancherel type estimates}\label{section2}

%The formulation of our general result requires the introduction of some notation. 
 Let $\om$ be a compactly supported smooth function defined on $\rn$.  
For  $\la\in\bbn_0$ let $\{\om^{\la}_{k}\}_{k\in\mathbb Z^n}$ be a sequence of 
compactly supported and smooth functions, defined on $\mathbb R^n$ by the formula $\om^\la_{k}(\xi)=2^{\la n/2}\om(2^\la \xi-k)$. These have the following properties:
\begin{enumerate}
\item[(i)] $ \{\om_{k}^{\la}\}_{k\in \mathbb Z^n}$ have  almost  disjoint  supports. 
\item[(ii)]  $\sum_{k\in\mathbb Z^n }{|\om_{k}^{\la}(\xi)|}\leq 2^{\lambda n/2}$ for all $\xi \in \mathbb R^n$.
\end{enumerate}
As a consequence of (i) and (ii) we obtain 
\begin{equation}\label{lrcondition}
\Big(\sum_{k\in\mathbb Z^n}\big|\om_{k}^{\la}(\xi)\big|^q\Big)^{1/q} \approx_q \sum_{k\in\mathbb Z^n}\big|\om_{k}^{\la}(\xi)\big|\leq 2^{\lambda n/2} 
\end{equation} 
for any $0<q<\infty$, due to the property of the supports.
We define
\begin{equation}\label{producttype}
\om_{\kkk}^{\la}(\xxxi):=\om_{k_1}^{\la}(\xi_1)\cdots\om_{k_m}^{\la}(\xi_m)
\end{equation} 
where $\kkk:=(k_1,\dots,k_m)\in (\bbzn)^m$ and $\xxxi=(\xi_1,\dots,\xi_m)\in(\bbrn)^m$.
Let $\UU$ be a subset of $(\bbzn)^m$ and $\{b^{\la}_{\kkk}\}_{\kkk\in \UU}$ be a sequence of complex numbers. 
We define 
\begin{equation}\label{sigmamultiplier}
\si^{\la}(\xxxi):=\sum_{\kkk\in\mathcal{U}} b^{\la}_{\kkk} \om_{\kkk}^{\la}(\xxxi)
\end{equation} and the corresponding $m$-linear multiplier operator by
\begin{equation}\label{multiplieroperator}
T_{\si^{\la}}\big(f_1,\dots,f_m\big)(x):=\int_{(\bbrn)^m}{\si^{\la}(\xxxi)\Big(\prod_{j=1}^{m}\wh{f_j}(\xi_j)\Big)e^{2\pi i\langle x,\sum_{j=1}^{m}\xi_j\rangle}}d\xxxi
\end{equation}
for $   x\in\bbrn$ and Schwartz functions $f_1,\dots,f_m$ on $\bbrn$. This operator 
coincides with that in \eqref{FTofW} when $\si^\la=\wh{W}$.

The multiplier $\si^\la$  defined in \eqref{sigmamultiplier} appears naturally in the decomposition of many  operators and plays a key role in the understanding of their boundedness. Actually in the bilinear case, one has 
  the estimate  (\cite{BGHH, Gr_He_Sl})
\begin{equation}\label{e10171}
\|T_{\si^\la}(f,g)\|_{L^1(\rn)}\le C\|\{b_{\vec k}^\la\}\|_{\ell^\nf}2^{\la n}| \mathcal U|^{1/4}\|f\|_{L^2(\rn)}\|g\|_{L^2(\rn)}.
\end{equation}
The presence of   $\|\{b_{\vec k}^\la\}\|_{\ell^\nf}$ and  
 $|U|^{1/4}$ indicate  the contribution of both the height and the support 
of $\si^\la$. This phenomenon is dissimilar to the $L^2$ boundedness of linear  multipliers, where the support of the multiplier plays no role. Motivated  by   many applications in which $\si^\la$ is an important building block, in this work we  obtain the $m$-linear version of \eqref{e10171}.

\begin{prop}\label{maintheorem}
Let $N$ be a positive integer and $\UU$ be a subset of $(\bbzn)^m$ with $|\UU|\le N$. For $\la\ge 0$, let $\{b^{\la}_k\}_{\kkk\in (\bbzn)^m}$ be a sequence of complex numbers satisfying $\Vert \{b^{\la}_{\kkk}\}_{\kkk}\Vert_{\ell^{\infty}}\le A_{\la}$. Let $\si^{\la}$ be defined as in (\ref{sigmamultiplier}).
Then there exists a constant $C= C_{n,m}>0$ such that
\begin{equation*}
\big\Vert T_{\si^{\la}}\big(f_1,\dots,f_m \big)\big\Vert_{L^{2/m}(\bbrn)}\le C A_{\la} N^{\frac{m-1}{2m}}2^{\frac{\la mn}{2}} \prod_{j=1}^{m}\Vert f_j\Vert_{L^2(\bbrn)}
\end{equation*}
for Schwartz functions $f_1,\dots,f_m$ on $\bbrn$.
\end{prop}

When $m=1$, Proposition~\ref{maintheorem} follows from     Plancherel's identity
and yields a bound on the $L^2$ norm of the corresponding linear operator $T_\si$
that depends only on the  height of the multiplier $\si$.  
For $m=2$, it coincides with \eqref{e10171}.
Below we focus on the consequences of Proposition~\ref{maintheorem} while its proof is postponed until the next section.

\begin{rmk}
After completing this paper, we were informed that Kato, Miyachi, and Tomita \cite{Katoetal1} recently obtained a result  that implies Proposition~\ref{maintheorem}. Their  proof is independent of ours and builds on their previous work in  \cite{Katoetal}.
\end{rmk}

%\vspace{0.2in}

The restriction $|\UU|\le N$ in Proposition~\ref{maintheorem} can be interpreted in terms of the compact support condition of $\si^{\la}$. Indeed, the support of $  \si^{\la} $ 
has measure bounded by a constant times $N$.
As we have seen in the proof of Proposition~\ref{maintheorem}, the $L^2\times\cdots\times L^2\to L^{2/m}$ boundedness of $m$-linear multiplier operator $T_{\si}$, $m\ge 2$, may be affected by the support of $\si$ while the $L^2$ boundedness depends only on $\Vert \si\Vert_{L^{\infty}}$ in the linear setting. 

On the other hand, the following ``support-independent" result could be obtained from Proposition~\ref{maintheorem} under an extra $\ell^q$ condition which is satisfied in many applications.

\begin{prop}\label{04231}
Let $m\in \bbn$ and  $0<q<\frac{2m}{m-1}$. Fixing $\la\in\mathbb N_0$, let $\{\om_{\kkk}^{\la}\}_{\kkk\in (\bbzn)^m}$ be wavelets of level $\la$.
Suppose $\{b^{\la}_{\kkk}\}_{\kkk\in (\bbzn)^m}$ is a sequence of complex numbers satisfying $\Vert \{b^{\la}_{\kkk}\}_{\kkk\in (\bbzn)^m}\Vert_{\ell^{\infty}}\leq A_{\la}$ and $\Vert \{b^{\la}_{\kkk}\}_{\kkk\in (\bbzn)^m}\Vert_{\ell^{q}}\leq B_{\la,q}$.
Then the $m$-linear multiplier $\sigma^{\la}$, defined in (\ref{sigmamultiplier}) with $\UU=(\bbzn)^m$, satisfies
\begin{equation*}
\big\Vert T_{\sigma^{\la}}(f_1,\dots,f_m)\big\Vert_{L^{2/m}}\lesssim A_{\la}^{1-\frac{(m-1)q}{2m}}B_{\la,q}^{\frac{(m-1)q}{2m}}2^{\la mn/2} \prod_{j=1}^{m}\Vert f_j\Vert_{L^2}
\end{equation*}
for Schwartz functions $f_1,\dots,f_m$ on $\bbrn$.
\end{prop}

\bpf
When $m=1$, it is clear from Plancherel's identity. Therefore we assume $m\ge 2$.

For $r\in \bbn$ let
$$
\UU^{\la}_r:=\big\{\kkk\in(\bbzn)^m: A2^{-r}<|b^{\la}_{\kkk}|\le A2^{-r+1}\big\}.
$$
As $\Vert \{b^{\la}_{\kkk}\}_{\kkk\in (\bbzn)^m}\Vert_{\ell^{\infty}}\le A_{\la}$,
  $(\bbzn)^m$ can be written as the disjoint union of $\UU^{\la}_r$, $r\in\bbn$,  
and thus we may decompose $\sigma^{\la}$ as 
$$\sigma^{\la}=\sum_{r\in\bbn}\sigma_r^{\la}$$
 where $\sigma_r^{\la}:=\sum_{\kkk\in \UU^{\la}_r}{b^{\la}_{\kkk}\om^{\la}_{\kkk}}$.
Observe that
\begin{equation}\label{bqrestrict}
2^{-r}A_{\la}|\UU^{\la}_r|^{1/q}\le \Big( \sum_{\kkk\in \UU^{\la}_r}{|b^{\la}_{\kkk}|^q}\Big)^{1/q}\le B_{\la,q},
\end{equation}
which implies 
\begin{equation}\label{cardur}
|\UU^{\la}_r|\le \Big( \frac{B_{\la,q}}{2^{-r}A_{\la}} \Big)^q.
\end{equation}
Applying Proposition~\ref{maintheorem} and (\ref{cardur}) to each $\sigma_r^{\la}$, we obtain
\begin{align*}
\big\Vert T_{\sigma_r^{\la}}(f_1,\dots,f_m)\big\Vert_{L^{2/m}}&\lesssim |\UU^{\la}_r|^{(m-1 )/2m}2^{\la mn/2}(A_{\la}2^{-r}) \prod_{j=1}^{m}\Vert f_j\Vert_{L^2}\\
 &\le (A_{\la} 2^{-r})^{1-\frac{(m-1)q}{2m}}B_{\la,q}^{\frac{q(m-1)}{2m}}2^{\la mn/2}\prod_{j=1}^{m}\Vert f_j\Vert_{L^2}.
\end{align*}
Taking $\ell^{2/m}$-norm over $r\in \bbn$, we have
\begin{align*}
\big\Vert T_{\sigma^{\la}}(f_1,\dots,f_m)\big\Vert_{L^{2/m}}&\le \Big( \sum_{r\in\bbn}{\big\Vert T_{\sigma_r^{\la}}(f_1,\dots,f_m)\big\Vert_{L^{2/m}}^{2/m}}\Big)^{m/2}\\
 &\lesssim A_{\la}^{1-\frac{(m-1)q}{2m}}B_{\la,q}^{\frac{(m-1)q}{2m}}2^{\la mn/2} \Vert f_1\Vert_{L^2}\cdots\Vert f_m\Vert_{L^2}
\end{align*} 
since $1-\frac{(m-1)q}{2m}>0$ and $2/m\le 1$.
\epf

For the case $m\ge 2$ and $q\ge \tf{2m}{m-1}$, we have the following substitute under an extra condition that all $\kkk$ belong to a ball of radius $C2^{\la}$, centered at the origin, which means that $\si^\la$ is contained in a ball of radius comparable to $1$.
%This result is useful in dealing with the multilinear multipliers of H\"ormander's type.

\begin{prop}\label{05071}
Let $m$ be a positive integer with $m\ge 2$ and $\frac{2m}{m-1}\le q< \infty$.
For each $\la\in \bbn_0$ let $\{\om_{\kkk}^{\la}\}$ be wavelets of level $\la$. 
Let $\UU^{\la}:=\{\kkk\in (\bbzn)^m: |\kkk|\leq C 2^{\la}\}$ for some $C>0$.
Suppose $\{b^{\la}_{\kkk}\}_{\kkk\in (\bbzn)^m}$ is a sequence of complex numbers with  $\Vert \{b^{\la}_{\kkk}\}_{\kkk\in (\bbzn)^m}\Vert_{\ell^{q}}\leq B_{\la,q}$.
Then the $m$-linear multiplier $\sigma^{\la}$, defined in (\ref{sigmamultiplier}) with $\UU=\UU^{\la}$, satisfies
\begin{equation*}
\big\Vert T_{\sigma^{\la}}(f_1,\dots,f_m)\big\Vert_{L^{2/m}}\lesssim  B_{\la,q} D_{\la,q,m}\Big(\prod_{j=1}^{m}\Vert f_j\Vert_{L^2}\Big)
\end{equation*}
for Schwartz functions $f_1,\dots,f_m$ on $\bbrn$,
where $$D_{\la,q,m}:=\begin{cases}
\la^{m/2}2^{\la mn/2}, & q=\frac{2m}{m-1}\\
2^{\la n(\frac{2m-1}{2}-\frac{m}{q})}, & q>\frac{2m}{m-1}
\end{cases}.$$
\end{prop}

\bpf
Pick $r_{max}\in \bbn$ satisfying $\frac{\la mn}{q} \le r_{max}<\frac{\la mn}{q}+1$.
Define
\begin{equation*}
\UU^{\la}_{r_{max}}:= \{\kkk\in\UU^{\la}: |b^{\la}_{\kkk}|\leq 2^{-r_{max}+1}B_{\la,q}\}
\end{equation*}
and for $1\le r< r_{max}$
\begin{equation*}
\UU_{r}^{\la}:=\{\kkk\in \UU^{\la}: 2^{-r}B_{\la,q}<|b^{\la}_{\kkk}|\le 2^{-r+1}B_{\la,q} \}.
\end{equation*}
Since $|b^{\la}_{\kkk}|\le B_{\la,q}$ for all $\kkk\in(\bbzn)^m$, we can write
\begin{equation*}
\si^{\la}=\sum_{r=1}^{r_{max}}\si_r^{\la}
\end{equation*} where $\si_r^{\la}:=\sum_{\kkk\in \UU_r^{\la}}{b^{\la}_{\kkk}\om_{\kkk}^{\la}}$.
Using the same argument in (\ref{bqrestrict}), we see that
\begin{equation}\label{cardurla1}
|\UU_r^{\la}|\le 2^{rq}, \qq 1\le r<r_{max}
\end{equation} and 
\begin{equation}\label{cardurla2}
|\UU_{r_{max}}^{\la}|\le |\UU^{\la}|\lesssim 2^{\la mn}\le 2^{r_{max}q}.
\end{equation}
Applying Proposition~\ref{maintheorem}, (\ref{cardurla1}), and (\ref{cardurla2}) to each $\sigma_{r}^{\la}$, we obtain
\begin{align*}
\big\Vert T_{\sigma_r^{\la}}(f_1,\dots,f_m)\big\Vert_{L^{2/m}}&\lesssim |\UU_r^{\la}|^{(m-1 )/2m}2^{\la mn/2}\Vert \{b^{\la}_{\kkk}\}_{\kkk\in \UU_r^{\la}}\Vert_{\ell^{\infty}}\prod_{j=1}^{m}\Vert f_j\Vert_{L^2}\\
 &\lesssim 2^{rq(m-1)/2m}2^{\la mn/2}2^{-r}B_{\la,q} \prod_{j=1}^{m}\Vert f_j\Vert_{L^2}\\
 &=B_{\la,q} 2^{\la mn/2}2^{r(\frac{q(m-1)}{2m}-1)}\prod_{j=1}^{m}\Vert f_j\Vert_{L^2}.
\end{align*}
Taking the $\ell^{2/m}$ norm over $1\le r\le r_{max}$, we have
\begin{align*}
\Vert T_{\sigma^{\la}}(f_1,\dots,f_m)\Vert_{L^{2/m}}&\le \Big( \sum_{r=1}^{r_{max}}{\big\Vert T_{\sigma_r^{\la}}(f_1,\dots,f_m)\big\Vert_{L^{2/m}}^{2/m}}\Big)^{m/2}\\
 &\lesssim B_{\la,q}2^{\la mn/2} \Big(\sum_{r=1}^{r_{max}}{2^{r(\frac{q(m-1)}{2m}-1)\frac{2}{m}}} \Big)^{m/2} \prod_{j=1}^{m}\Vert f_j\Vert_{L^2}.
\end{align*}
We note that
\begin{align*}
\Big(\sum_{r=1}^{r_{max}}{2^{r(\frac{q(m-1)}{2m}-1)\frac{2}{m}}} \Big)^{m/2}&\approx \begin{cases}
r_{max}^{m/2}, & q=\frac{2m}{m-1}\\
  2^{r_{max}(\frac{q(m-1)}{2m}-1)},        & q>\frac{2m}{m-1}
\end{cases}\\
&\lesssim \begin{cases}
\la^{m/2}, & q=\frac{2m}{m-1}\\
  2^{\la n(\frac{m-1}{2}-\frac{m}{q})},     & q>\frac{2m}{m-1}
\end{cases},
\end{align*}
which completes the proof.
\epf

In the study of bilinear rough singular integrals and bilinear H\"ormander multipliers, an argument splitting the problem to diagonal and off-diagonal cases is utilized. The off-diagonal case uses  Plancherel's theorem and a pointwise control. In the diagonal case, we employ the bilinear result of Plancherel type, which is actually the driving force of this work.
We now present a multilinear version generalizing and combining these two parts, which shows that all $l$-linear Plancherel type result,  $1\le l\le m$, is necessary in the study of many $m$-linear multipliers.

For $\mu\in\bbn_0$ let $\VV^{\mu}$ be a subset of $\{\kkk\in (\bbzn)^m:2^{\mu-c_0}\le |\kkk|\le 2^{\mu+c_0} \}$ for some $c_0\ge 1$.
Let $M$ be a positive constant and for each $1\le l\le m$ let
\begin{equation*}
\VV_{l}^{\mu}:=\big\{ \kkk\in \VV^{\mu}:  |k_1|,\dots,|k_l|\ge M>|k_{l+1}|,\dots,|k_m|\big\}.
\end{equation*} 
%In other words, $\VV_{l}^{\mu}$ consists of $M^{m-l}$ number of $l$-columns. 
%If $\kkk\in \VV_{l}^{\mu+\la}$, then $\om_{\kkk}$ is supported in an annulus of radius $2^{\mu}$.
We also define  $L^{\la}_kf:=\big( \om^{\la}_k \wh{f}\,\big)^{\vee}$ and
$L_k^{\la,\ga}f:=\big(\om^{\la}_k(\cdot/2^{\ga})\wh{f}\, \big)^{\vee}$ for $k\in\bbzn$.

\begin{prop}\label{keyapplication1}
Let $m$ be a positive integer with $m\ge 2$ and $0<q<\frac{2m}{m-1}$.
For each $\la\in\bbn_0$, let $\{\om_{\kkk}^{\la}\}_{\kkk}$ be wavelets of level $\la$.
Suppose that $\{b_{\kkk}^{\la,\ga,\mu}\}_{\ga,\mu\in \bbz, \kkk\in (\bbzn)^m}$ is a sequence of complex numbers satisfying 
\begin{equation*}
\sup_{\ga \in \bbz}\big\Vert \{b_{\kkk}^{\la,\ga,\mu}\}_{\kkk\in (\bbzn)^m} \big\Vert_{\ell^{\infty}}\le A_{\la,\mu}
\end{equation*}
 and 
 \begin{equation*}
 \sup_{\ga \in\bbz}\big\Vert \{b_{\kkk}^{\la,\ga,\mu}\}_{\kkk\in (\bbzn)^m} \big\Vert_{\ell^{q}}\le B_{\la,\mu,q}.
 \end{equation*}
Then the following statements hold:
 \begin{enumerate}
 \item  For $1\le r\le 2$ there exists a constant $C>0$, independent of $\la, \mu$, such that
\begin{align*}
&\Big\Vert \Big(\sum_{\ga\in\bbz}\Big| \sum_{\kkk\in \VV_1^{\la+\mu}} b_{\kkk}^{\la,\ga,\mu}L_{k_1}^{\la,\ga}f_1^{\la,\ga,\mu}\prod
_{j=2}^{m}L_{k_j}^{\la,\ga}f_j\Big|^r\Big)^{1/r}\Big\Vert_{L^{2/m}}\\
&\le C A_{\la,\mu} 2^{\la mn/2}\Big(\sum_{\ga\in\bbz}{\Vert f_1^{\la,\ga,\mu}\Vert_{L^2}^r} \Big)^{1/r} \prod_{i=2}^{m}\Vert f_{i}\Vert_{L^2}
\end{align*}  
for Schwartz functions $f_1,\dots,f_m$ on $\bbrn$.

 \item For $2\le l\le m$ there exists a constant $C>0$, independent of $\la, \mu$, such that
\begin{align*}
&\Big\Vert \sum_{\ga\in\bbz} \Big| \sum_{\kkk\in \VV_l^{\la+\mu}} b_{\kkk}^{\la,\ga,\mu}\Big( \prod_{j=1}^{l}L_{k_j}^{\la,\ga}f_j^{\la,\ga,\mu}\Big)\Big( \prod_{j=l+1}^{m}L_{k_{j}}^{\la,\ga}f_{j}\Big) \Big| \Big\Vert_{L^{2/m}}\\
&\le C A_{\la,\mu}^{1-\frac{(l-1)q}{2l}}B_{\la,\mu,q}^{\frac{(l-1)q}{2l}} 2^{\la mn/2}\Big[ \prod_{j=1}^{l}\Big(\sum_{\ga\in\bbz}{\Vert f_j^{\la,\ga,\mu}\Vert_{L^2}^2} \Big)^{1/2}\Big] \Big[\prod_{j=l+1}^{m}\Vert f_{j}\Vert_{L^2}\Big]
\end{align*} 
for Schwartz functions $f_1,\dots,f_m$ on $\bbrn$,
 where $\prod_{m+1}^m$ is understood as empty.
 \end{enumerate}
\end{prop}

The proof of Proposition~\ref{keyapplication1} is given in the next section following   that of Proposition~\ref{maintheorem}.

%\section{Applications of the main theorem}\label{09031}

\section{Proofs of Proposition~\ref{maintheorem} and Proposition~\ref{keyapplication1}}\label{pfmain}

When $m=1$, Proposition~\ref{maintheorem} follows immediately from Plancherel's identity. Thus, we will be concerned only with the case $m\ge 2$.
For the bilinear case $m=2$, a concept called column is used; see, for instance, \cite{BGHH, Gr_He_Sl}. A column $Col_k^{\mathcal U}$ with $k\in\bbz$ related to a   subset $\mathcal U\sset \bbz^2$ is defined as $\mathcal U\cap \{\{k\}\times\bbz\}$. We  generalize this concept to higher dimensions in least two ways, expressed in terms of the 
dimension and the codimension of the set. A column related to $\mathcal U\sset\bbz^3$ could be $\mathcal U\cap \{(k,l)\times\bbz\}$ of dimension $1$, or $\mathcal U\cap \{\{k\}\times\bbz^2\}$ of dimension $2$. It turns out that   both notions are needed to handle the case $m=3$.
%, and we call the first tyep $1$-column, and the second type $2$-column in terms of their dimensions.

We introduce several notions and study their combinatorial properties.
For a  fixed $\kkk\in (\bbzn)^m$, $1\leq l\leq m$, and $1\leq j_1\le  \dots\le j_l\le m$ let 
$$\kkk^{j_1,\dots,j_l}:=(k_{j_1},\dots,k_{j_l})$$ denote the vector in $(\bbzn)^l$ consisting of $j_1,\dots,j_l$ components of $\kkk$ and 
$\kkk^{*j_1,j_2,\dots,j_l}$ stand for the vector in $(\bbzn)^{m-l}$, consisting of $\kkk$ excepting $j_1$, \dots, $j_l$ components
(e.g.  $\kkk^{*1,\dots,j}=\kkk^{j+1,\dots,m}=(k_{j+1},\dots,k_m)\in (\bbzn)^{m-j}$).
For any sets $\UU$ in $(\bbzn)^m$, $1\le j\le m$, and $1\le j_1\le\dots\le j_l \leq m$ let 
$$\mathcal{P}_j\UU:=\big\{k_j\in\bbzn:\kkk \in\UU ~\text{ for some }\kkk^{*j}\in(\bbzn)^{m-1} \big\}$$
$$\mathcal{P}_{*j_1,\dots,j_l}\UU:=\big\{\kkk^{*j_1,\dots,j_l}\in(\bbzn)^{m-l}:\kkk\in \UU ~\text{ for some }k_{j_1},\dots, k_{j_l}\in\bbzn \big\}$$ 
be the projections of $\UU$ onto the $k_j$-column and $\kkk^{*j_1,\dots,j_l}$-plane, respectively.
For a fixed $\kkk^{*j_1,\dots,j_l}\in \mathcal{P}_{*j_1,\dots,j_l}\UU$,  we define 
\begin{equation*}
Col_{\kkk^{*j_1,\dots,j_l}}^{\UU}:=\{\kkk^{j_1,\dots,j_l}\in (\bbzn)^l: \kkk=(k_1,\dots,k_m) \in\UU\}.
\end{equation*}
Then we observe that 
\begin{equation}\label{disjointdecom1}
\sum_{\kkk\in \UU}\cdots = \sum_{\kkk^{*j_1,\dots,j_l}\in \mathcal{P}_{*j_1,\dots,j_l}\UU}\Big( \sum_{\kkk^{j_1,\dots,j_l}\in Col_{\kkk^{*j_1,\dots,j_l}}^{\UU}}\cdots\Big).
\end{equation}
Furthermore, for each $\kkk^{*j_1,\dots,j_{l-1}}$ we have
$$Col_{\kkk^{*j_1,\dots,j_l}}^{\UU}=\bigcup_{k_{j_l}\in \mathcal{P}_{j_l}Col_{\kkk^{*j_1,\dots,j_l}}^{\UU}}Col_{\kkk^{*j_1,\dots,j_{l-1}}}^{\UU}\times \{k_{j_l}\}$$
and this allows us to wrtie
\begin{equation}\label{disjointdecom2}
\sum_{\kkk\in \UU}\cdots =\!\!\!\sum_{\kkk^{*j_1,\dots,j_l}\in \mathcal{P}_{*j_1,\dots,j_l}\UU}\Big( \sum_{k_{j_l} \in \mathcal{P}_{j_l}Col_{\kkk^{*j_1,\dots,j_l}}^{\UU}}{\Big(\sum_{\kkk^{j_1,\dots,j_{l-1}}\in Col_{\kkk^{*j_1,\dots,j_{l-1}}}^{\UU}}\!\!\!\!\!\!
\!\!\!\!\cdots  \Big)}\Big).
\end{equation}

To make it easier to understand, let us think about 
%For example, we consider
 the case $n=1$ and $m=3$. $\mathcal P_1\mathcal U$ is the projection of $\mathcal U$ to the first coordinate. $\mathcal P_{*1}\mathcal U$ is the projection of $\mathcal U$ to the $(k_2,k_3)$-plane. $Col_{\kkk^{*1}}^{\mathcal U}$ is an $1$-column in $\bbz^3\cap \mathcal U$ with $(k_2,k_3)={\kkk}^{*1}$ fixed. When $j_l=1$, 
identity \eqref{disjointdecom1} says that 
$$
\sum_{\kkk \in\mathcal U}\cdots =\sum_{(k_2,k_3)\in\mathcal P_{*1}\mathcal U}\q \Big(\sum_{k_1\in Col^{\mathcal U}_{k_2,k_3}} \cdots \Big).
$$
When $(j_1,\dots, j_l)=(1,2)$,
identity \eqref{disjointdecom2} says that
$$
\sum_{\kkk \in\mathcal U}\cdots =\sum_{k_3\in\mathcal P_3\mathcal U}\q \Big( \sum_{k_2\in \mathcal P_2Col^{\mathcal U}_{k_3}}\Big( \sum_{k_1\in Col^{\mathcal U}_{k_2,k_3}}  \cdots   \Big) \Big).
$$

The proof of Proposition~\ref{maintheorem} is based on the   decompositions in (\ref{disjointdecom1}) and (\ref{disjointdecom2})  and on the following lemma.

\begin{lm}\label{keylemma}
Let $m\ge 2$ and $\mathcal{U}$ be a subset of $(\bbzn)^m$. Let $\la\in\bbn_0$ and $\{\om_{\kkk}^{\la}\}_{\kkk\in(\bbzn)^m}$ be wavelets whose level is $\la$. 
Let $\si^{\la}=\sum_{\vec k\in\mathcal U}b^{\la}_{\vec k}\om^{\la}_{\vec k}$, where $\{b^{\la}_{\kkk}\}_{\kkk\in\mathcal{U}}$ is a sequence of complex numbers satisfying $\Vert \{ b^{\la}_{\kkk} \}_{\kkk\in (\bbzn)^m}\Vert_{\ell^{\infty}} \le A_{\la}$.
Then there exists a constant $C>0$ such that
\begin{align*}
\big\|T_{\si^{\la}}(f_1,\cdots,f_m)\big\|_{L^{2/m}}&\le C A_{\la}2^{\la (m-1) n/2} \Big(\prod_{i\ne j,1\le i\le m}{\Vert f_i\Vert_{L^2}}\Big)\\
 &\qq \times \Big( \int_{\bbrn}{\big|\wh{f_j}(\xi)\big|^2\sum_{\kkk\in \mathcal{U}}{\big|\om^{\la}_{k_j}(\xi)\big|^2}}d\xi\Big)^{1/2}
\end{align*}
for each $1\le j\le m$.

\end{lm}

 \bpf
Without loss of generality, we may assume $j=1$. 
In view of (\ref{disjointdecom1}), $\sigma^{\la}$ can be written as
\begin{equation*}
\sigma^{\la}(\xxxi)=\sum_{\kkk^{*1}\in \mathcal{P}_{*1}\mathcal{U}}\om^{\la}_{k_2}(\xi_2)\cdots\om^{\la}_{k_m}(\xi_m)\sum_{k_1\in Col^{\mathcal{U}}_{\kkk^{*1}}}{b^{\la}_{\kkk}\om^{\la}_{k_1}(\xi_1)}, 
\end{equation*} 
and this yields that
\begin{equation*}
T_{\sigma^{\la}}\big(f_1,\dots,f_m\big)(x)=\sum_{\kkk^{*1}\in \mathcal{P}_{*1}\mathcal{U}}\Big(\prod_{i=2}^{m}L^{\la}_{k_i}f_i(x)\Big)\sum_{k_1\in Col^{\mathcal{U}}_{\kkk^{*1}}}{b^{\la}_{\kkk}L^{\la}_{k_1}f_1(x)},
\end{equation*}
where $L^{\la}_{k}f:=(\om^{\la}_k \wh{f}\, )^{\vee}$ for $k\in\bbzn$.
Using the Cauchy-Schwarz inequality  and H\"older's inequality, we obtain that
\begin{align*}
\big\Vert T_{\sigma^{\la}}(f_1,\dots,f_m)\big\Vert_{L^{2/m}}&\le \Big\Vert \Big(\sum_{\kkk^{*1}\in (\bbzn)^{m-1}}\Big| \prod_{i=2}^{m}L^{\la}_{k_i}f_i \Big|^2 \Big)^{1/2}\Big\Vert_{L^{2/(m-1)}}\\
 &\qq \times \Big\Vert \Big( \sum_{\kkk^{*1}\in \mathcal{P}_{*1}\mathcal{U}}{\Big| \sum_{k_1\in Col^{\mathcal{U}}_{\kkk^{*1}}}{b^{\la}_{\kkk}L^{\la}_{k_1}f_1}\Big|^2}\Big)^{1/2}\Big\Vert_{L^2}\\
 &=:I\times II.
\end{align*}
As a direct consequence of Plancherel's identity and (\ref{lrcondition}), we have
\begin{equation*}
\big\Vert \big\{L^{\la}_kf \big\}_{k\in\bbzn}\big\Vert_{L^2(\ell^2)}\lesssim 2^{\la n/2}\Vert f\Vert_{L^2}
\end{equation*}
and thus,
\begin{align*}
I&=\Big\Vert \prod_{i=2}^{m}\Big(\sum_{k_i\in\bbzn}{|L^{\la}_{k_i}f_i|^2} \Big)^{1/2}\Big\Vert_{L^{2/(m-1)}} \leq \prod_{i=2}^{m}{\big\Vert \big\{ L^{\la}_{k_i}f_i\big\}_{k_i\in\bbzn}\big\Vert_{L^2(\ell^2)}}\\
&\lesssim 2^{\la (m-1)n/2}\prod_{i=2}^{m}\Vert f_i\Vert_{L^2}, 
\end{align*} 
where  the first inequality is obtained by H\"older's inequality.
Moreover, it follows from Plancherel's identity and the disjoint compact support property of $\{\om^{\la}_{k_1}\}_{k_1\in\bbzn}$ that
\begin{align*}
II&\lesssim \Big(\sum_{\kkk^{*1}\in \mathcal{P}_{*1}\mathcal{U}}{\Big\Vert \wh{f_1}\sum_{k_1\in Col^{\mathcal{U}}_{\kkk^{*1}}}b^{\la}_{\kkk}\om^{\la}_{k_1}\Big\Vert_{L^2}^2} \Big)^{1/2}\\
          &\approx \Big( \int_{\bbrn}{|\wh{f_1}(\xi)|^2\sum_{\kkk^{*1}\in\mathcal{P}_{*1}\mathcal{U}}\sum_{k_1\in Col^{\mathcal{U}}_{\kkk^{*1}}}{|b^{\la}_{\kkk}|^2|\om^{\la}_{k_1}(\xi)|^2} }d\xi\Big)^{1/2}
\end{align*}
and this is controlled by a constant multiple of
\begin{equation*}
A_{\la}\Big( \int_{\bbrn}{|\wh{f_1}(\xi)|^2\sum_{\kkk\in \mathcal{U}}{|\om^{\la}_{k_1}(\xi)|^2}}d\xi\Big)^{1/2}
\end{equation*} where (\ref{disjointdecom1}) is applied.
This completes the proof.
\epf

%We now prove Proposition~\ref{maintheorem}. 
\subsection{Proof of Proposition~\ref{maintheorem}}
Let $N_1,\dots,N_{m-1}$ be positive numbers less than $N$, which will be chosen later.
We separate $\UU$ into $m$ disjoint subsets
\begin{align*}
\UU^1&:=\{\kkk\in\mathcal{U}:|Col_{\kkk^{*1}}^{\mathcal{U}}|> N_1\}\\
\UU^2&:=\{\kkk\in\mathcal{U}\setminus \mathcal{U}^1:|Col_{\kkk^{*1,2}}^{\UU}|> N_2\}\\
&\vdots\\
\UU^{m-1}&:=\{\kkk\in\mathcal{U}\setminus (\mathcal{U}^1\cup\cdots\cup\mathcal{U}^{m-2}):|Col_{\kkk^{*1,\dots,m-1}}^{\UU}|> N_{m-1}\}\\
\UU^m&:=\mathcal{U}\setminus (\mathcal{U}^1\cup\cdots\cup\mathcal{U}^{m-1})
\end{align*}
and write
\begin{equation*}
\sigma^{\la}=\sum_{j=1}^{m}\sum_{\kkk\in\mathcal{U}^{j}}{b^{\la}_{\kkk}\om^{\la}_{\kkk}}=:\sum_{j=1}^{m}\sigma^{\la}_{(j)}.
\end{equation*}

Observe that for $1\le j\le m-1$, due to (\ref{disjointdecom1}),
$$N\ge |\UU^j|> N_j|\mathcal{P}_{*1,\dots,j}\UU^j|,$$ 
which implies 
\begin{equation}\label{cardp}
|\mathcal{P}_{*1,\dots,j}\UU^j|< NN_j^{-1}.
\end{equation}
Moreover, for $2\le j\le m$ and $\kkk^{*1,\dots,j-1}\in \mathcal{P}_{*1,\dots,j-1}\UU^j$,
\begin{equation}\label{cardc}
|Col_{\kkk^{*1,\dots,j-1}}^{\UU^j}|\leq N_{j-1},
\end{equation}  which follows from the fact $\UU^j\subset \UU\setminus \UU^{j-1}$.

We now apply Lemma \ref{keylemma} to each $\sigma_{(j)}$, $1\leq j\leq m$, to obtain 
\begin{align}\label{lemmaapply}
\big\Vert T_{\si^{\la}_{(j)}}(f_1,\cdots,f_m)\big\Vert_{L^{2/m}}&\le C A_{\la} 2^{\la (m-1) n/2} \prod_{i\ne j,1\le i\le m}{\Vert f_i\Vert_{L^2}}\\
 &\qq \times \Big( \int_{\bbrn}{|\wh{f_j}(\xi)|^2\sum_{\kkk\in \mathcal{U}^j}{|\om^{\la}_{k_j}(\xi)|^2}}d\xi\Big)^{1/2}\nonumber.
\end{align}
Note that
\begin{equation*}
\sum_{\kkk\in \UU^1}{|\om^{\la}_{k_1}(\xi)|^2}=\sum_{\kkk^{*1}\in\mathcal{P}_{*1}\UU^1}{\Big(\sum_{k_1\in Col_{\kkk^{*1}}^{\UU^1}}{|\om^{\la}_{k_1}(\xi)|^2} \Big)}\leq 2^{\la n}|\mathcal{P}_{*1}\UU^1|<2^{\la n}NN_1^{-1}
\end{equation*} where (\ref{disjointdecom1}), (\ref{lrcondition}), and (\ref{cardp}) are applied.
Similarly, when $2\le j\le m-1$,  we have
\begin{align*}
\sum_{\kkk\in \UU^j}|\om^{\la}_{k_j}(\xi)|^2&=\sum_{\kkk^{*1,\dots,j}\in\mathcal{P}_{*1,\dots,j}\UU^j}{\Big( \sum_{k_j\in\mathcal{P}_jCol_{\kkk^{*1,\dots,j}}^{\UU^j}}{|\om^{\la}_{k_j}(\xi)|^2|Col_{\kkk^{*1,\dots,j-1}}^{\UU^j}|}\Big)}\\
&\le 2^{\la n}N_{j-1}|\mathcal{P}_{*1,\dots,j}\UU^j| \le 2^{\la n} NN_{j-1}N_{j}^{-1},
\end{align*} using (\ref{disjointdecom2}), (\ref{lrcondition}), (\ref{cardc}), and (\ref{cardp}).
For the last case $j=m$, it follows from (\ref{disjointdecom1}), (\ref{lrcondition}), and (\ref{cardc}) that
\begin{equation*}
\sum_{\kkk\in\UU^m}|\om^{\la}_{k_m}(\xi)|^2=\sum_{k_m\in\mathcal{P}_m\UU^m}{|\om^{\la}_{k_m}(\xi)|^2|Col_{k_m}^{\UU^m}|}\leq 2^{\la n}N_{m-1}.
\end{equation*}

Now we choose $N_1,\dots,N_{m-1}$ satisfying the identity
\begin{equation}\label{nidentity}
NN_1^{-1}=NN_1N_2^{-1}=NN_2N_3^{-1}=\cdots=
NN_{m-2}N_{m-1}^{-1}
=N_{m-1}.
\end{equation}
Solving (\ref{nidentity}), we have 
$$N_j=N^{j/m}, \qq 1\le j\le m-1$$
and this establishes
$$\sum_{\kkk\in\UU^j}{|\om^{\la}_{k_j}(\xi)|^2}\leq 2^{\la n}N^{(m-1)/m}, \qq 1\leq j\leq m,$$
which further implies
\begin{equation*}
\Big( \int_{\bbrn}{|\wh{f_j}(\xi)|^2\sum_{\kkk\in \mathcal{U}^j}{|\om^{\la}_{k_j}(\xi)|^2}}d\xi\Big)^{1/2}\le 2^{\la n/2}N^{(m-1)/2m}\Vert f_j\Vert_{L^2}.
\end{equation*}

Then this, together with (\ref{lemmaapply}), proves
\begin{equation*}
\big\Vert T_{\sigma^{\la}_{(j)}}(f_1,\dots,f_m)\big\Vert_{L^{2/m}}\lesssim A_{\la}N^{(m-1)/2m}2^{\la mn/2}\prod_{i=1}^{m}\Vert f_i\Vert_{L^2}
\end{equation*}
as desired.

\subsection{Proof of Proposition~\ref{keyapplication1}}

%\section{Main techniques in multilinear operators}\label{basictechniques}

We first observe that
\begin{equation}\label{cardpv}
\big| \mathcal{P}_{*1,\dots,l}\VV_{l}^{\mu}\big|\le M^{n(m-l)} \qq \text{for }~ \mu\ge 0,
\end{equation}
\begin{equation*}
L_k^{\la,\ga}f(x)=L^{\la}_k\big(f(\cdot/2^{\ga})\big)(2^{\ga}x),
\end{equation*}
and
\begin{equation*}
|L_k^{\la,\ga}f(x)|\lesssim 2^{\la n/2}\mathcal{M}f(x)\qquad \text{ for }~k\in\bbzn.
\end{equation*}
Here $\mathcal{M}$ is the Hardy-Littlewood maximal operator, defined by $\mathcal{M}f(x):= \sup_{Q\ni x} |Q|^{-1} \int_{Q}{|f(y)|}dy$,  where the supremum is taken over all cubes containing $x$.
Then in view of (\ref{disjointdecom1})  we can write
\begin{align}\label{mainptest}
&\hspace{.3in}\Big| \sum_{\kkk\in\VV_l^{\la+\mu}}{b_{\kkk}^{\la,\ga,\mu} \Big(\prod_{j=1}^{l}   L_{k_j}^{\la,\ga}f_j^{\la,\ga,\mu}(x)\Big) \Big( \prod_{j=l+1}^{m}L_{k_{j}}^{\la,\ga}f_{j}(x)\Big)  }\Big| \nonumber\\
&\hspace{.2in}\lesssim 2^{\la (m-l)n/2} \Big[\prod_{j=l+1}^{m}\mathcal{M}f_{j}(x)\Big]\\
&\q  \times \sum_{\kkk^{*1,\dots,l}\in \mathcal{P}_{*1,\dots,l}\VV_{l}^{\la+\mu}}\Big|\sum_{\kkk^{1,\dots,l}\in Col_{\kkk^{*1,\dots,l}}^{\VV_{l}^{\la+\mu}}}b_{\kkk}^{\la,\ga,\mu}\Big(\prod_{j=1}^{l}L^{\la}_{k_j}\big( f_j^{\la,\ga,\mu}(\cdot/2^{\ga})\big)(2^{\ga}x)\Big)\Big|\nonumber.
\end{align}

When $l=1$, using (\ref{mainptest}), H\"older's inequality, (Minkowski inequality for $r<2$ ), the $L^2$ boundedness of $\mathcal{M}$, and (\ref{cardpv}), we obtain
\begin{align*}
&\Big\Vert \Big(\sum_{\ga\in\bbz}\Big| \sum_{\kkk\in \VV_1^{\la+\mu}} b_{\kkk}^{\la,\ga,\mu}L_{k_1}^{\la,\ga}f_1^{\la,\ga,\mu} \Big( \prod_{j=2}^{m}L_{k_j}^{\la,\ga}f_j\Big)\Big|^r\Big)^{1/r}\Big\Vert_{L^{2/m}}\\
&\lesssim 2^{\la (m-1)n/2}\Big(\prod_{j=2}^{m}\big\Vert f_j\big\Vert_{L^2} \Big) \\
&\qq \times \sum_{\kkk^{*1}\in \mathcal{P}_{*1}\VV_{1}^{\la+\mu}}\Big( \sum_{\ga\in\bbz}\Big\Vert   \sum_{k_1\in Col_{\kkk^{*1}}^{\VV_{1}^{\la+\mu}}}b_{\kkk}^{\la,\ga,\mu}L^{\la}_{k_1}\big(f_1^{\la,\ga,\mu}(\cdot/2^{\ga}) \big)(2^{\ga}\cdot)      \Big\Vert_{L^2}^r\Big)^{1/r}.
\end{align*}
A change of variables and Plancherel's identity yield that
\begin{equation*}
\Big\Vert   \sum_{k_1\in Col_{\kkk^{*1}}^{\VV_{1}^{\la+\mu}}}b_{\kkk}^{\la,\ga,\mu}L^{\la}_{k_1}\big( f_1^{\la,\ga,\mu}(\cdot/2^{\ga}) \big)(2^{\ga}\cdot)      \Big\Vert_{L^2}\le A_{\la,\mu}2^{\la n/2}\Vert f_1^{\la,\ga,\mu}\Vert_{L^2},
\end{equation*}
which proves the first estimate.

Similarly, for $2\le l\le m$,  we can see
\begin{align*}
&\Big\Vert \sum_{\ga\in\bbz}  \Big| \sum_{\kkk\in \VV_l^{\la+\mu}} b_{\kkk}^{\la,\ga,\mu}\Big( \prod_{j=1}^{l}L_{k_j}^{\la,\ga}f_j^{\la,\ga,\mu}\Big) \Big( \prod_{j=l+1}^{m}L_{k_{j}}^{\la,\ga}f_{j}\Big) \Big|\Big\Vert_{L^{2/m}}\\
&\lesssim 2^{\la (m-l)n/2}\Big( \prod_{j=l+1}^{m}\Vert f_{j}\Vert_{L^2}\Big)\q\times  \\
& \sum_{\kkk^{*1,\dots,l}\in \mathcal{P}_{*1,\dots,l}\VV_{l}^{\la+\mu}} \! 
 \Big\Vert \sum_{\ga\in\bbz}  \Big|\!\!\sum_{\kkk^{1,\dots,l}\in Col_{\kkk^{*1,\dots,l}}^{\VV_{l}^{\la+\mu}}} \!\!\! b_{\kkk}^{\la,\ga,\mu}\Big[ \prod_{j=1}^{l}L^{\la}_{k_i}\big( f_j^{\la,\ga,\mu}(\cdot/2^{\ga}) \big)(2^{\ga}\cdot)\Big]\Big|       \Big\Vert_{L^{2/l}}.
\end{align*}
The $L^{2/l}$ norm is clearly dominated by 
\begin{equation*}
\Big(\sum_{\ga\in\bbz}\Big\Vert   \sum_{\kkk^{1,\dots,l}\in Col_{\kkk^{*1,\dots,l}}^{\VV_{l}^{\la+\mu}}}b_{\kkk}^{\la,\ga,\mu} \Big( \prod_{j=1}^{l}L^{\la}_{k_j}\big( f_j^{\la,\ga,\mu}(\cdot/2^{\ga}) \big)(2^{\ga}\cdot)\Big)        \Big\Vert_{L^{2/l}}^{2/l} \Big)^{l/2}
\end{equation*} since $2/l\le 1$,
and now we apply a change of variables, Proposition \ref{04231}, and H\"older's inequality to obtain that the above expression is less than
\begin{align*}
&A_{\la,\mu}^{1-\frac{(l-1)q}{2l}}B_{\la,\mu,q}^{\frac{(l-1)q}{2l}}2^{\la ln/2}\Big( \sum_{\ga\in\bbz}2^{-\ga n}\Big(\prod_{j=1}^{l}\big\Vert f_j^{\la,\ga,\mu}(\cdot/2^{\ga})\big\Vert_{L^2}^{2/l}\Big) \Big)^{l/2} \\
&\le A_{\la,\mu}^{1-\frac{(l-1)q}{2l}}B_{\la,\mu,q}^{\frac{(l-1)q}{2l}}2^{\la ln/2}\prod_{j=1}^{l}\Big( \sum_{\ga\in\bbz}{\Vert f_j^{\la,\ga,\mu}\Vert_{L^2}^2}\Big)^{1/2}.
\end{align*} This completes the proof.

\section{Compactly supported wavelets}

Typical functions possessing  properties (i) and (ii) in Section~\ref{section2} 
are the compactly supported wavelets constructed by Daubechies \cite{Dau1}; their construction is contained in the books of Meyer \cite{Me} and Daubechies \cite{Dau2}.
Wavelets have been used to study singular integrals in different settings; see for instance \cite{Meyer3}, \cite{Hyt}, \cite{DPWW}, and  \cite{Gr_He_Ho}.
For the purposes of this paper, we need smooth wavelets with compact supports but also of product type,  like \eqref{producttype}.
%; this means that they have the form \eqref{producttype} described above. 
The construction of such 
orthonormal bases  
is carefully presented in Triebel \cite{Tr2010}, but for the reader's sake we provide an outline.
For any fixed $M\in\bbn$ there exist real compactly supported functions $\psi_F, \psi_M$ in  $\mathscr{C}^M(\bbr)$ satisfying the following properties:
\begin{enumerate}
\item[(a)] $\Vert \psi_F\Vert_{L^2(\bbr)}=\Vert \psi_M\Vert_{L^2(\bbr)}=1$
\item[(b)] $\int_{\bbr}{x^{\alpha}\psi_M(x)}dx=0$ for all $0\le \alpha \le M$
\item[(c)] If      $\Psi_{\GGG}$ is a function on $\bbr^{mn}$, defined by
$$\Psi_{\GGG}(\xxx):=\psi_{g_1}(x_1)\cdots \psi_{g_{mn}}(x_{mn})$$
for $\xxx:=(x_1,\dots,x_{mn})\in \bbr^{mn}$ and $\GGG:=(g_1,\dots,g_{mn})$ in the set     
 $$\II:=\big\{\GGG:=(g_1,\dots,g_{mn}):g_i\in\{F,M\} \big\},$$
then the family of functions
\begin{equation*}
\bigcup_{\la\in\bbn_0}\bigcup_{\kkk\in \bbz^{mn}}\big\{ 2^{\la{mn/2}}\Psi_{\GGG}(2^{\la}\xxx-\kkk):\GGG\in \II^{\la}\big\}
\end{equation*}
forms an orthonormal basis of $L^2(\bbr^{mn})$,
where $\II^0:=\II$ and $\II^{\la}:=\II\setminus \{(F,\dots,F)\}$ for $\la \ge 1$.
\end{enumerate}

Fix  $1<q<\infty$ and $s\ge 0$. Let  $\| F \|_{ L_s^q(\bbr^{mn})}$ denote the Sobolev space norm defined as the $ L^q((\bbrn)^m)$ norm 
of $(\vec{I}-\vec{\Delta})^{s/2}F$, where $\vec{\Delta}$ is the Laplacian of a function $F$ on $(\bbrn)^m$.
It is also shown in \cite{Tr2006} that
if $M$ is sufficiently large and $F$ is a tempered distribution on $\bbr^{mn}$ lying in $L^q_s(\mathbb R^{mn})$, then $F$ can be represented as
\begin{equation}\label{daubechdecomp}
F(\xxx)=\sum_{\la\in\bbn_0}\sum_{\GGG\in\II^{\la}}\sum_{\kkk\in \bbz^{mn}}b_{\GGG,\kkk}^{\la}2^{\la mn/2}\Psi_{\GGG}(2^{\la} \xxx -\kkk)
\end{equation}
and 
$$
\Big\Vert \Big(\sum_{\vec G,\ \vec k} \big|b_{\vec G,\vec k}^\la \Psi^\la_{\vec G,\vec k}\big|^2\Big)^{1/2} \Big\Vert_{L^q(\bbr^{mn})} \le C2^{-s\la}\|F\|_{L^q_s(\bbr^{mn})} , 
$$
where $ \Psi^\la_{\vec G,\vec k}(\vec x)=2^{\la mn/2}\Psi_{\GGG}(2^{\la} \xxx -\kkk)$, and 
\begin{equation*}
b_{\GGG,\kkk}^{\la}:=\int_{\bbr^{mn}}{F(\xxx)\Psi^\la_{\vec G,\vec k}(\vec x)}d\xxx.
\end{equation*}
Moreover, it follows from the last estimate and disjoint supports of $\Psi^\la_{\vec G,\vec k}$ that 
\begin{align}
\big\Vert \big\{b_{\GGG,\kkk}^{\la}\big\}_{\kkk\in \bbz^{mn}}\big\Vert_{\ell^{q}}
\approx&\Big(2^{\la mn(1-q/2)}\int_{\bbr^{mn}}\Big( \sum_{\vec k} \big|b^\la_{\vec G,\vec k}\Psi^\la_{\vec G,\vec k}(\vec x)\big|^2\Big)^{q/2}d\vec x\Big)^{1/q}\notag \\
\lesssim& \; 2^{-\la (s-mn/q+mn/2)}\Vert F\Vert_{L^q_s(\bbr^{mn})}. \label{lqestimate}
\end{align} 

\medskip

%\noindent{\bf Notation:} 
We will write %$\xxxi:=(\xi_1,\dots,\xi_m)\in (\bbrn)^m$, 
$\GGG:=(G_1,\dots,G_m)\in (\{F,M\}^n)^m$, and
\begin{equation*}
\Psi_{\GGG}(\xxxi)=\Psi_{G_1}(\xi_1)\cdots \Psi_{G_m}(\xi_m).
\end{equation*}
For each $\kkk:=(k_1,\dots,k_m)\in (\bbzn)^m$ and $\la\in \bbn_0$, let 
$$
\Psi_{G_i,k_i}^{\la}(\xi_i):=2^{\la n/2}\Psi_{G_i}(2^{\la}\xi_i-k_i), \qq 1\le i\le m
$$
and
$$
\Psi_{\GGG,\kkk}^{\la}(\xxxi):=\Psi_{G_1,k_1}^{\la}(\xi_1)\cdots \Psi_{G_m,k_m}^{\la}(\xi_m).
$$
We also assume that the  support of $\psi_{g_i}$ is contained in $\{\xi\in \bbr: |\xi|\le C_0 \}$ for some $C_0>1$,
which implies that
\begin{equation}\label{supppsi}
\textup{Supp}(\Psi_{G_i,k_i}^\la)\subset \big\{\xi_i\in\bbrn: |2^{\la}\xi_i-k_i|\le C_0\sqrt{n}\big\}.
\end{equation}
In other words, the support of $\Psi_{G_i,k_i}^\la$ is contained in the ball centered at $2^{-\la}k_i$ and radius $C_0\sqrt{n}2^{-\la}$.

\section{Proof of Theorem \ref{application1}}\label{pfapp1}

Using (\ref{daubechdecomp}) with $s=0$, we decompose $\si$ as
\begin{equation*}
\si(\xxxi)=\sum_{\la\in\bbn_0}\sum_{\GGG\in\II^{\la}}\sum_{\kkk\in (\bbzn)^m}{b_{\GGG,\kkk}^{\la}\Psi_{G_1,k_1}^{\la}(\xi_1)\cdots \Psi_{G_m,k_m}^{\la}(\xi_m)}=:\sum_{\la\in\bbn_0}\sum_{\GGG\in \II^{\la}}\si_{\GGG}^{\la}(\xxxi)
\end{equation*}
where $b_{\GGG,\kkk}^{\la}:= \int_{(\bbrn)^m}{\si(\xxxi)\Psi_{\GGG,\kkk}^{\la}(\xxxi)}d\xxxi $.
As an immediate consequence of Proposition \ref{04231}, we have
\begin{equation*}
\big\Vert T_{\si_{\GGG}^{\la}}(f_1,\dots,f_m)\big\Vert_{L^{2/m}}\lesssim \big\Vert \{b_{\GGG,\kkk}^{\la}\}_{\kkk}\big\Vert_{\ell^{\infty}}^{1-\frac{(m-1)q}{2m}} \big\Vert \{b_{\GGG,\kkk}^{\la}\}_{\kkk}\big\Vert_{\ell^{q}}^{\frac{(m-1)q}{2m}}2^{\la mn/2} \prod_{j=1}^{m}\Vert f_j\Vert_{L^2}.
\end{equation*}
We first observe that (\ref{lqestimate}) yields that
\begin{equation*}
\big\Vert \{b_{\GGG,\kkk}^{\la}\}_{\kkk}\big\Vert_{\ell^{q}}\lesssim 2^{\la mn(1/q-1/2)}\Vert \si\Vert_{L^q((\bbrn)^m)}.
\end{equation*}
In addition, as $\si\in \mathscr{C}^{M_q}((\bbrn)^{m})$, using this property, the $M_q$ vanishing moment condition of $\Psi_{\GGG,\kkk}^{\la}$ in conjunction with Taylor's formula, an argument similar to \cite[Lemma 2.1]{Gr_He_Sl} and  \cite[Lemma 7]{Gr_He_Ho} yields 
\begin{equation*}
\big\Vert \{b_{\GGG,\kkk}^{\la}\}_{\kkk}\big\Vert_{\ell^{\infty}}\lesssim 2^{-\la (M_q+mn/2)}D_0.
\end{equation*}
Therefore, we finally arrive at the estimate
\begin{align*}
&\big\Vert T_{\si_{\GGG}^{\la}}(f_1,\dots,f_m)\big\Vert_{L^{2/m}}\\
&\lesssim 2^{-\la \big(M_q(1-\frac{(m-1)q}{2m})-\frac{ n(m-1)}{2}\big)}D_0^{1-\frac{(m-1)q}{2m}}\Vert \si\Vert_{L^q((\bbrn)^m)}^{\frac{(m-1)q}{2m}} \prod_{j=1}^{m}\Vert f_j\Vert_{L^2},
\end{align*}
which in turn implies that
\begin{align*}
&\big\Vert T_{\si}(f_1,\dots,f_m)\big\Vert_{L^{2/m}}\\
&\le \Big( \sum_{\la\in\bbn_0}\sum_{\GGG\in\II^{\la}}{\big\Vert T_{\si_{\GGG}^{\la}}(f_1,\dots,f_m)\big\Vert_{L^{2/m}}^{2/m}}\Big)^{m/2}\\
&\lesssim D_0^{1-\frac{(m-1)q}{2m}}\Vert \si\Vert_{L^q((\bbrn)^m)}^{\frac{(m-1)q}{2m}}  \Big( \sum_{\la\in\bbn_0}{2^{-\la (M_q(1-\frac{(m-1)q}{2m})-\frac{ n(m-1)}{2})\frac{2}{m}}}\Big)^{m/2}\prod_{j=1}^{m}\Vert f_j\Vert_{L^2}.
\end{align*}
Since $M_q>\frac{ m(m-1)n}{2m-(m-1)q}$, the sum over $\la$   converges and 
completes the proof.

\section{Proof of Theorem \ref{application2}}\label{pfapp2}

Without loss of generality, we may assume $\frac{2m}{m+1}<q<2$ as $L^r(\mathbb S^{mn-1})\sset L^q(\mathbb S^{mn-1})$ for $r\ge q$.
We first utilize a dyadic decomposition introduced by Duoandikoetxea and Rubio de Francia \cite{Du_Ru}.
Recall that $\Phi^{(m)}$ is a Schwartz function such that $\wh{\Phi^{(m)}}$ is supported in the annulus $\{\xxxi\in (\bbrn)^m: 1/2\le |\xxxi|\le 2\}$ and $\sum_{j\in\bbz}\wh{\Phi^{(m)}_j}(\xxxi)=1$ for $\xxxi\not= \0$
where $\wh{\Phi_j^{(m)}}(\xxxi):=\wh{\Phi^{(m)}}(\xxxi/2^j)$.

 For $\ga\in\bbz$  let
 $$ K^{\ga}(\yyy):=\wh{\Phi^{(m)}}(2^{\ga}\yyy)K(\yyy), \quad \yyy\in (\bbrn)^m$$
 and then we observe that $K^\ga(\yyy)=2^{\ga mn} K^0(2^\ga \yyy)$.
 For $\mu\in\bbz$ we define
 $$
K_{\mu}^{\ga}(y):=\Phi_{{\mu}+\ga}^{(m)}\ast K^{\ga}(y)=2^{\ga mn}[\Phi_{{\mu}}^{(m)}\ast K^{0}](2^\ga y).$$
It follows from this definition that 
$$
\wh{K^\ga_\mu}(\vec \xi)= \wh{\Phi^{(m)}}(2^{-(\mu+\ga)}\vec \xi)\wh{K^0}(2^{-\ga}\vec \xi)=\wh{K^0_\mu}(2^{-\ga}\vec \xi),
$$
which implies that $\wh{K^\ga_\mu}$ is bounded uniformly in $\ga$ while they have almost disjoint supports, so it is natural to add them together as follows,
 $$K_{\mu}(\yyy):=\sum_{\ga\in \bbz}{K_{\mu}^{\ga}(\yyy)}.$$

 We define
 \begin{equation*}
 \LL_{\mu}\big(f_1,\dots,f_m\big)(x):=\int_{(\bbrn)^m}{K_{\mu}(\yyy)\prod_{j=1}^{m}f_j(x-y_j)} ~ d\yyy, \q x\in\bbrn
 \end{equation*} and write
 \begin{align}
 \big\Vert \LL_{\Om}(f_1,\dots,f_m)\big\Vert_{L^{2/m}}&\lesssim \Big\Vert \sum_{{\mu}\in\bbz : {2^{{\mu}-10}\le C_0 \sqrt{mn}}}{\LL_{\mu}(f_1,\dots,f_m)}\Big\Vert_{L^{2/m}}\nonumber\\
  &\qq +\Big\Vert \sum_{{\mu}\in\bbz :2^{{\mu}-10}> C_0\sqrt{mn}}{\LL_{\mu}(f_1,\dots,f_m)}\Big\Vert_{L^{2/m}} \label{secondterm}
 \end{align} where $C_0$ is the constant that appeared in (\ref{supppsi}).

 First of all, using the argument in \cite[Proposition 3]{Gr_He_Ho}, we can prove that 
 \begin{equation}\label{e09022}
 \big\Vert \LL_{\mu}(f_1,\dots,f_m)\big\Vert_{L^p}\lesssim \Vert \Omega\Vert_{L^q(\mathbb{S}^{mn-1})}\Big( \prod_{j=1}^{m}\Vert f_j\Vert_{L^{p_j}}\Big) \begin{cases}
2^{(mn-\delta){\mu}}, & {\mu}\ge 0\\
2^{(1-\delta){\mu}}, & {\mu}<0
 \end{cases}
 \end{equation}  
 for $0<\delta<1/q'$, and this implies that
 \begin{equation*}
 \Big\Vert \sum_{{\mu\in\bbz}:2^{{\mu}-10}\le C_0\sqrt{mn}}{\LL_{\mu}(f_1,\dots,f_m)}\Big\Vert_{L^{2/m}}\lesssim \Vert \Omega\Vert_{L^q(\mathbb{S}^{mn-1})}\prod_{j=1}^{m}\Vert f_j\Vert_{L^2}.
 \end{equation*}
 It remains to bound the term (\ref{secondterm}), but this can be reduced to proving that for $2^{\mu-10}>C_0\sqrt{mn}$, there exists $\epsilon_0>0$ such that
 \begin{equation}\label{mainmaingoal}
 \big\Vert {\LL_{\mu}(f_1,\dots,f_m)}\big\Vert_{L^{2/m}}\lesssim 2^{-\epsilon_0 \mu}\Vert \Om \Vert_{L^q(\mathbb{S}^{mn-1})}\prod_{j=1}^{m}\Vert f_j\Vert_{L^2},
 \end{equation}  which compensate the estimate (\ref{e09022}) for $\mu \ge 0$.
 Recall that
 \begin{equation*}
 \wh{K_{\mu}}(\xxxi)=\sum_{\ga\in\bbz}{\wh{K_{\mu}^0}(\xxxi/2^{\ga})}
 \end{equation*}
 and 
\begin{equation}\label{cptsupportk}
\textup{Supp}\wh{K_{\mu}^{0}}\subset \big\{\xxxi\in (\bbrn)^m: 2^{{\mu}-1}\le |\xxxi|\le 2^{{\mu}+1} \big\}.
\end{equation}
Using (\ref{daubechdecomp}), $\wh{K_{\mu}^0}$ can be written as
\begin{equation}\label{kj0}
\wh{K_{\mu}^0}(\xxxi)=\sum_{\la\in\bbn_0}\sum_{\GGG\in\II^{\la}}\sum_{\kkk\in (\bbzn)^m}b_{\GGG,\kkk}^{\la,\mu}\Psi_{G_1,k_1}^{\la}(\xi_1)\cdots \Psi_{G_m,k_m}^{\la}(\xi_m)
\end{equation}
where 
\begin{equation*}
b_{\GGG,\kkk}^{\la,\mu}:=\int_{(\bbrn)^m}{\wh{K_{\mu}^0}(\xxxi)\Psi_{\GGG,\kkk}^{\la}(\xxxi)}d\xxxi.
\end{equation*}

It is already known in \cite[Lemma 7]{Gr_He_Ho} that
\begin{equation}\label{maininftyest}
\big\Vert \{b_{\GGG,\kkk}^{\la,\mu}\}_{\kkk}\big\Vert_{\ell^{\infty}}\lesssim      2^{-\delta {\mu}}2^{-\la (M+1+mn)} \Vert \Om\Vert_{L^q(\mathbb{S}^{mn-1})}
\end{equation} where $M$ is the number of vanishing moments of $\Psi_{\GGG}$ and $0<\delta<1/q'$.
In addition, (\ref{lqestimate}), the Hausdorff-Young inequality, and Young's inequality prove that
\begin{align}\label{mainlqest}
\big\Vert \{b_{\GGG,\kkk}^{\la,\mu}\}_{\kkk}\big\Vert_{\ell^{q'}}&\lesssim 2^{-\la mn (1/2-1/q')}\Vert \wh{K_{\mu}^0}\Vert_{L^{q'}}\nonumber\\
&\lesssim 2^{-\la mn(1/q-1/2)}\Vert \Omega \Vert_{L^q(\mathbb{S}^{mn-1})}.
\end{align}

Furthermore, if $2^{{\mu}-10}>C_0\sqrt{mn}$, then we may replace $\kkk\in (\bbzn)^m$ in (\ref{kj0}) by $2^{\la+{\mu}-2}\le |\kkk|\le 2^{\la+{\mu}+2}$
due to (\ref{cptsupportk}) and the compact support condition of $\Psi_{\GGG,\kkk}^{\la}$.
Therefore the proof of (\ref{mainmaingoal}) can be reduced to the inequality
\begin{align}\label{mmainest}
&\Big\Vert  \sum_{\la\in\bbn_0}\sum_{\GGG\in\II^{\la}}\sum_{\ga\in\bbz}\sum_{\kkk\in\UU^{\la+{\mu}}}b_{\GGG,\kkk}^{\la,\mu}\prod_{j=1}^{m}L_{G_j,k_j}^{\la,\ga}f_j     \Big\Vert_{L^{2/m}}\nonumber\\
&\lesssim 2^{-\epsilon_0\mu}\Vert \Om\Vert_{L^q(\mathbb{S}^{mn-1})}\prod_{j=1}^{m}\Vert f_j\Vert_{L^2}
\end{align}
where
\begin{equation*}
\UU^{\la+{\mu}}:=\big\{ \kkk\in (\bbzn)^m: 2^{\la+{\mu}-2}\le |\kkk|\le 2^{\la+{\mu}+2},~  |k_1|\ge\cdots \ge |k_m| \big\}
\end{equation*}
and 
\begin{equation}\label{lgklg}
L_{G,k}^{\la,\ga}f:=\big( \Psi_{G,k}^{\la}(\cdot/2^{\ga})\wh{f}\big)^{\vee}.
\end{equation}
 Here, it is additionally assumed that $|k_1|\ge \dots\ge |k_m|$ in $\UU^{\la+\mu}$ as the remaining  cases follow by symmetry. 
%We remark that 
%\begin{equation}\label{e09021}
%\sum_{\la\in\bbn_0}\sum_{\GGG\in\II^{\la}}\sum_{\ga\in\bbz}\sum_{\kkk\in\UU^{\la+{\mu}}}b_{\GGG,\kkk}^{{\mu},\la}\prod_{j=1}^{m}L_{G_j,k_j}^{\la,\ga}f_j  \in L^{2/m}(\rn)
%\end{equation}
%by \eqref{e09022} since the operator is essentially $\mathcal L_\mu(f_1,\dots, f_m)$.
Then we note that $\UU^{\la+{\mu}}$ can be expressed as the union of $m$ disjoint subsets
\begin{equation*}
\UU^{\la+{\mu}}_1:=\{\kkk\in \UU^{\la+{\mu}}:|k_1|\ge 2C_0\sqrt{n}>|k_2|\ge\cdots\ge |k_m|      \}
\end{equation*}
\begin{equation*}
\UU^{\la+{\mu}}_2:=\{\kkk\in \UU^{\la+{\mu}}:|k_1|\ge |k_2|\ge 2C_0\sqrt{n}>|k_3|\ge\cdots\ge |k_m|      \}
\end{equation*}
$$\vdots$$
\begin{equation*}
\UU^{\la+{\mu}}_m:=\{\kkk\in \UU^{\la+{\mu}}:|k_1|\ge \cdots\ge |k_m|\ge 2C_0\sqrt{n}      \}.
\end{equation*}

The function in the left-hand side of  (\ref{mmainest})   could be written as
\begin{equation*}
 \sum_{l=1}^{m}  \sum_{\la\in\bbn_0}\sum_{\GGG\in\II^{\la}}  \sum_{\ga\in\bbz}  \TT_{\GGG,l}^{\la,\ga,\mu}\big(f_1,\dots,f_m\big)        
  \end{equation*}
where
\begin{equation}\label{toperator}
\TT_{\GGG,l}^{\la,\ga,\mu}\big(f_1,\dots,f_m\big):=\sum_{\kkk\in\UU_l^{\la+{\mu}}}b_{\GGG,\kkk}^{\la,\mu} \Big(\prod_{j=1}^{m}L_{G_j,k_j}^{\la,\ga}f_j \Big).
\end{equation}
Observe that when $\kkk\in \UU_l^{\la+{\mu}}$, 
\begin{equation}\label{chaequi}
L_{G_j,k_j}^{\la,\ga}f_j=L_{G_j,k_j}^{\la,\ga}f_j^{\la,\ga,{\mu}} \qq \text{for }~ 1\le j\le l
\end{equation}
due to the support of $\Psi_{G_j,k_j}^{\la}$,
where $\wh{f_j^{\la,\ga, {\mu}}}(\xi_j):=\wh{f_j}(\xi_j)\chi_{C_0\sqrt{n}2^{\ga-\la}\le |\xi|\le 2^{\ga+{\mu}+3}}$.
Moreover, it is easy to show that   for $\mu\ge 10$ and $\la\in\bbn_0$, 
\begin{equation}\label{L2}
\Big( \sum_{\ga\in\bbz}{\big\Vert f_j^{\la,\ga,{\mu}}\big\Vert_{L^2}^2}\Big)^{1/2}\lesssim ({\mu}+\la)^{1/2}\Vert f_j\Vert_{L^2} \lesssim \mu^{1/2}(\la+1)^{1/2}\Vert f_j\Vert_{L^2} 
\end{equation} 
where Plancherel's identity is applied   in the first inequality.

Now we claim that for each $1\le l\le m$ there exists $\epsilon_0, M_0>0$ such that
\begin{align}
&\Big\Vert  \sum_{\ga\in\bbz}\TT_{\GGG,l}^{\la,\ga,\mu}\big(f_1,\dots,f_m\big) \Big\Vert_{L^{2/m}}\lesssim  \Vert \Omega\Vert_{L^q(\mathbb{S}^{mn-1})}2^{- \epsilon_0\mu}2^{-\la M_0}\prod_{j=1}^{m}\Vert f_j\Vert_{L^2}.\label{finalgoal}
\end{align} 
Then the left-hand side of (\ref{mmainest}) is controlled by a constant times
\begin{align*}
&\Big(  \sum_{l=1}^{m}  \sum_{\la\in\bbn_0}\sum_{\GGG\in\II^{\la}}  \Big\Vert  \sum_{\ga\in\bbz}\TT_{\GGG,l}^{\la,\ga,\mu}\big(f_1,\dots,f_m\big) \Big\Vert_{L^{2/m}}^{2/m}    \Big)^{m/2}\\
&\lesssim 2^{-\epsilon_0\mu}   \Vert \Omega\Vert_{L^q(\mathbb{S}^{mn-1})}\prod_{j=1}^{m}\Vert f_j\Vert_{L^2},
\end{align*}
which completes the proof of (\ref{mmainest}). Therefore, it remains to prove (\ref{finalgoal}).

\subsection{ The case $l=1$} 
The proof relies on the fact that
 if $\widehat{g_{\ga}}$ is supported on $\{\xi\in \rn : C^{-1} 2^{\ga+\mu}\leq |\xi|\leq C2^{\ga+\mu}\}$ for some $C>1$ and $\mu\in \bbz$, then
\begin{equation}\label{marshallest}
\Big\Vert \Big\{ \Phi^{(1)}_j\!\ast\! \Big(\sum_{\ga\in\bbz}{g_{\ga}}\!\Big)\!\Big\}_{\! j\in\mathbb{Z}}\Big\Vert_{L^p(\ell^q)}\lesssim_{C} \big\Vert \big\{ g_j\big\}_{j\in\mathbb{Z}}\big\Vert_{L^p(\ell^q)} \q \text{uniformly in }~\mu
\end{equation}  for $0<p<\infty$. The proof of (\ref{marshallest}) is elementary and standard, so it is omitted here. See \cite[(13)]{Gr_He_Ho} and  \cite[Theorem 3.6]{Ya}  for a related argument.

Note that 
\begin{equation*}
2^{\la+{\mu}-3}\le 2^{\la+{\mu}-2}-2C_0\sqrt{mn}\le |\kkk|- (|k_2|^2+\dots+|k_m|^2)^{1/2}\le |k_1|\le 2^{\la+{\mu}+2}
\end{equation*}
and this implies that 
\begin{equation*}
\textup{Supp}\big(\Psi_{G_1,k_1}^{\la}(\cdot/2^{\ga}) \big)\subset \{\xi\in \bbrn: 2^{\ga+{\mu}-4}\le |\xi|\le 2^{\ga+{\mu}+3}\}.
\end{equation*}
Moreover,   since $|k_j|\le 2C_0\sqrt{n}$  for $2\le j\le m$ and $2^{{\mu}-10}>C_0\sqrt{mn}$,
\begin{equation*}
\textup{Supp}\big(\Psi_{G_j,k_j}^{\la}(\cdot/2^{\ga})\big)\subset \{\xi\in \bbrn: |\xi|\le m^{-1/2}2^{\ga+{\mu}-8}\}.
\end{equation*} 
Therefore, the Fourier transform of $\TT_{\GGG,1}^{\la,\ga,\mu}\big(f_1,\dots,f_m\big)$  is supported in the set $\{\xi\in \bbrn: 2^{\ga+{\mu}-5}\le |\xi|\le 2^{\ga+{\mu}+4} \}$.
Using the Littlewood-Paley theory for Hardy spaces \cite[Theorem 2.2.9]{MFA}, there exists a  unique  polynomial $Q^{\la,\mu,\GGG}(x)$ such that 
\begin{align}\label{polyargument}
&\Big\Vert  \sum_{\ga\in\bbz} \TT_{\GGG,1}^{\la,\ga,\mu}\big(f_1,\dots,f_m\big) -Q^{\la,\mu,\GGG}\Big\Vert_{L^{2/m}}\nonumber\\
&\lesssim \Big\Vert  \Big\{     \Phi^{(1)}_j \ast \Big( \sum_{\ga\in\bbz} \TT_{\GGG,1}^{\la,\ga,\mu}\big(f_1,\dots,f_m\big) \Big)\Big\}_{j\in\bbz} \Big\Vert_{L^{2/m}(\ell^2)}
\end{align}
and then (\ref{marshallest}) and (\ref{chaequi}) yield that the above $L^{2/m}(\ell^{2})$-norm is dominated by a constant multiple of 
\begin{equation*}
\Big\Vert \Big( \sum_{\ga\in\bbz}{\big| \TT_{\GGG,1}^{\la,\ga,\mu}\big(f_1^{\la,\ga,{\mu}},f_2,\dots,f_m\big)  \big|^2}\Big)^{1/2}\Big\Vert_{L^{2/m}}.
\end{equation*}
We now apply Proposition \ref{keyapplication1}, (\ref{maininftyest}), and (\ref{L2}) to bound the $L^{2/m}$-norm by 
\begin{align*}
&\big\Vert \big\{ b_{\GGG,\kkk}^{\la,\mu}\big\}_{\kkk}\big\Vert_{\ell^{\infty}}2^{\la mn/2}\Big( \sum_{\ga\in\bbz}{\big\Vert f_1^{\la,\ga,\mu}\big\Vert_{L^2}^2}\Big)^{1/2}\prod_{j=2}^{m}\Vert f_j\Vert_{L^2}\\
&\lesssim \Vert \Om\Vert_{L^q(\mathbb{S}^{mn-1})}2^{-\delta {\mu}}\mu^{1/2}2^{-\la(M+1+ mn/2)}(\la+1)^{1/2}\prod_{j=1}^{m}\Vert f_j\Vert_{L^2}
\end{align*}

This implies that the left-hand side of (\ref{polyargument}) is bounded by $$\Vert \Om\Vert_{L^q(\mathbb{S}^{mn-1})}2^{-\epsilon_0\mu}2^{-\la M_0}\prod_{j=1}^{m}\Vert f_j\Vert_{L^2}$$ 
for some $0<\epsilon_0<\delta$ and $0<M_0<M+1+mn/2$,
and thus  
\begin{equation}\label{polynomiall2m}
 \sum_{\ga\in\bbz} \TT_{\GGG,1}^{\la,\ga,\mu}\big(f_1,\dots,f_m\big) -Q^{\la,\mu,\GGG} \in L^{2/m}.
 \end{equation}

Furthermore, it follows from Proposition \ref{keyapplication1} that
\begin{align*}
& \Big\Vert  \sum_{\ga\in\bbz} \TT_{\GGG,1}^{\la,\ga,\mu}\big(f_1,\dots,f_m\big) \Big\Vert_{L^{2/m}}\le \Big\Vert  \sum_{\ga\in\bbz}\big|  \TT_{\GGG,1}^{\la,\ga,\mu}\big(f_1,\dots,f_m\big)\big| \Big\Vert_{L^{2/m}} \\
&\lesssim \big\Vert \big\{ b_{\GGG,\kkk}^{\la,\mu}\big\}_{\kkk}\big\Vert_{\ell^{\infty}} 2^{\la mn/2}\Big(\sum_{\ga\in \bbz}\Vert f_1^{\ga,\la,\mu}\Vert_{L^2} \Big)\prod_{j=1}^{m}\Vert f_j\Vert_{L^2}\\
&\lesssim 2^{-\de \mu} 2^{-\la(M+1+mn/2)}\Vert \Omega\Vert_{L^q(\mathbb{S}^{mn-1})}\Big(\sum_{\ga\in \bbz}\Vert f_1^{\la,\ga,\mu}\Vert_{L^2} \Big)\prod_{j=1}^{m}\Vert f_j\Vert_{L^2}.
\end{align*}
Since $f_1$ is a Schwartz function, we have
\begin{align}\label{l=1case}
\Big\Vert f_1^{\la,\ga,\mu}\Big\Vert_{L^2}&=\Big\Vert \wh{f_1^{\la,\ga,\mu}}\Big\Vert_{L^2}=\Big(\int_{C_0\sqrt{n} 2^{\ga-\la}\le |\xi|\le 2^{\ga+\mu+3}}{|\wh{f_1}(\xi)|^2}d\xi \Big)^{1/2}\nonumber\\
 &\lesssim_N \begin{cases}
 2^{(\ga+\mu)n/2}, ~& \ga <0\\
 2^{-(\ga-\la)(N-n/2)}, ~& \ga\ge 0
 \end{cases}
\end{align} for sufficiently large $N>n/2$, which yields that
$$
\sum_{\ga\in \bbz}\Vert f_1^{\la,\ga,\mu}\Vert_{L^2}
$$
is finite (of course, this depends on $\ga$, $\mu$, and $f_1$). Therefore, we also have
\begin{equation}\label{casel=1}
 \sum_{\ga\in\bbz} \TT_{\GGG,1}^{\la,\ga,\mu}\big(f_1,\dots,f_m\big) \in L^{2/m}
\end{equation}
and thus the polynomial $Q^{\la,\mu,\GGG}$ in (\ref{polynomiall2m}) should be zero.
In conclusion,
\begin{align*}
& \Big\Vert  \sum_{\ga\in\bbz} \TT_{\GGG,1}^{\la,\ga,\mu}\big(f_1,\dots,f_m\big) \Big\Vert_{L^{2/m}}\\
  &\lesssim \Vert \Om\Vert_{L^q(\mathbb{S}^{mn-1})}2^{-\epsilon_0 {\mu}}2^{-\la M_0}\prod_{j=1}^{m}\Vert f_j\Vert_{L^2}.
\end{align*}
This proves (\ref{finalgoal}) for $l=1$.

\subsection{The case $2\le l\le m$}
We apply (\ref{chaequi}), Proposition \ref{keyapplication1}  with $2<q'<\frac{2m}{m-1}$, (\ref{L2}), (\ref{maininftyest}), and (\ref{mainlqest}) to obtain that
\begin{align}
&\Big\Vert  \sum_{\ga\in\bbz}\TT_{\GGG,l}^{\la,\ga,\mu}\big(f_1,\dots,f_m\big) \Big\Vert_{L^{2/m}}\le \Big\Vert  \sum_{\ga\in\bbz} \big| \TT_{\GGG,l}^{\la,\ga,\mu}\big(f_1,\dots,f_m\big) \big| \Big\Vert_{L^{2/m}}\notag\\
&\lesssim \big\Vert \big\{ b_{\GGG,\kkk}^{\la,\mu}\big\}_{\kkk}\big\Vert_{\ell^{\infty}}^{1-\frac{(m-1)q'}{2m}} \big\Vert \big\{ b_{\GGG,\kkk}^{\la,\mu}\big\}_{\kkk}\big\Vert_{\ell^{q'}}^{\frac{(m-1)q'}{2m}} 2^{\la mn/2}(\la+1)^{l/2}{\mu}^{l/2}\prod_{j=1}^{m}\Vert f_j\Vert_{L^2}\notag\\
&\lesssim \Vert \Om\Vert_{L^q(\mathbb{S}^{mn-1})}2^{-\delta {\mu}(1-\frac{(m-1)q'}{2m})}{\mu}^{m/2}2^{-\la C_{M,m,n,q}}(\la+1)^{m/2}\prod_{j=1}^{m}\Vert f_j\Vert_{L^2}\label{e09024}
\end{align} 
where $$C_{M,m,n,q}:=(M+1+mn)(1-\frac{(m-1)q'}{2m})+mn(1/q-1/2)\frac{(m-1)q'}{2m}-\frac{mn}{2}.$$
Here we used the embedding $\ell^{q'}\hookrightarrow \ell^{\infty}$ and  the fact that $\frac{l-1}{2l}\le \frac{m-1}{2m}$.
  Then (\ref{finalgoal}) follows from choosing $M$ sufficiently large so that $C_{M,m,n,q}>0$ since $1-\frac{(m-1)q'}{2m}>0$.

\section{Proof of Theorem \ref{application3}}\label{pfapp3}

The strategy in this section is similar to that used in  handling multilinear rough singular integrals in Section \ref{pfapp2}, but the decomposition is more delicate. We describe the decomposition first.
Write
\begin{equation*}
\si(\xxxi)=\sum_{\ga\in\bbz}{\si_{\ga}(\xxxi/2^{\ga})}
\end{equation*}
where  $\si_{\ga}(\xxxi):=\sigma(2^{\ga} \xxxi)\wh{\Phi^{(m)}}(\xxxi)$.
Clearly,
 \begin{equation}\label{sigamma}
 \textup{Supp}(\si_{\ga})\subset \{\xxxi\in (\bbzn)^m: 1/2\le |\xxxi|\le 2\}
 \end{equation} and according to (\ref{daubechdecomp}), 
\begin{equation}\label{sigmagamma}
\si_{\ga}(\xxxi)=\sum_{\la\in \bbn_0}\sum_{\GGG\in \II^{\la}}\sum_{\kkk\in (\bbzn)^m}{b_{\GGG,\kkk}^{\la,\ga}\Psi_{G_1,k_1}^{\la}(\xi_1)\cdots \Psi_{G_m,k_m}^{\la}(\xi_m)}
\end{equation}
where $b_{\GGG,\kkk}^{\la,\ga}:=\int_{(\bbrn)^m}{\si_{\ga}(\xxxi)\Psi_{\GGG,\kkk}^{\la}(\xxxi)}d\xxxi$.
Moreover, it follows from (\ref{lqestimate}) that for $1<q<\infty$ and $s\ge 0$
\begin{equation}\label{lrestimate}
\big\Vert \{b_{\GGG,\kkk}^{\la,\ga}\}_{\kkk\in (\bbzn)^m}\big\Vert_{\ell^q} \lesssim   2^{-\la(s-mn/q+mn/2)} \big\Vert \sigma(2^{\ga}\;\vec{\cdot }\;)\wh{\Phi^{(m)}}\big\Vert_{L^q_s((\bbrn)^m)}  .
\end{equation}
As we did in the proof of Theorem \ref{application2}, it is enough to consider only the case $|k_1|\ge\dots\ge |k_m|$.
Therefore, we replace $\kkk\in (\bbzn)^m$ in (\ref{sigmagamma}) by $\kkk\in \UU:=\{\kkk\in (\bbzn)^{m}:|k_1|\ge \cdots\ge |k_m|\}$ and write
\begin{align*}
\si_{\ga}(\xxxi)&={\sum_{\la\in\bbn_0: 2^{\la}\ge 2^8 C_0m\sqrt{n}   } \;\sum_{\GGG\in \II^{\la}}\sum_{\kkk\in \UU}{b_{\GGG,\kkk}^{\la,\ga}\Psi_{G_1,k_1}^{\la}(\xi_1)\cdots \Psi_{G_m,k_m}^{\la}(\xi_m)} }\\
    &\qq +{\sum_{\la\in\bbn_0: 2^{\la}<2^8 C_0m\sqrt{n}   } \; \sum_{\GGG\in \II^{\la}}\sum_{\kkk\in \UU}{b_{\GGG,\kkk}^{\la,\ga}\Psi_{G_1,k_1}^{\la}(\xi_1)\cdots \Psi_{G_m,k_m}^{\la}(\xi_m)} }\\
    &=:\si_{\ga}^{(1)}(\xxxi)+\si_{\ga}^{(2)}(\xxxi).
\end{align*}
We are only concerned with $\si_{\ga}^{(1)}$ as a similar and simpler argument is applicable to the other one since the sum over $\la$ in $\si_{\ga}^{(2)}$ is finite sum.

If $2^8   C_0m\sqrt{n}    \le 2^{\la}$, then $b_{\GGG,\kkk}^{\la,\ga}$ vanishes unless $2^{\la-2}\le |\kkk|\le 2^{\la+2}$ due to (\ref{sigamma}) and the compact support of $\Psi_{\GGG}$.
Thus, letting
\begin{equation*}
\UU^{\la}:=\{\kkk\in\UU:2^{\la-2}\le |\kkk|\le 2^{\la+2}\},
\end{equation*}
we write
\begin{equation*}
\si_{\ga}^{(1)}(\xxxi)={\sum_{\la: 2^{\la}\ge 2^8C_0m\sqrt{n}} \; \sum_{\GGG\in \II^{\la}}\sum_{\kkk\in \UU^{\la}}{b_{\GGG,\kkk}^{\la,\ga}\Psi_{G_1,k_1}^{\la}(\xi_1)\cdots \Psi_{G_m,k_m}^{\la}(\xi_m)} }.
\end{equation*}

Now we split $\UU^{\la}$ into $m$ disjoint subsets
\begin{equation*}
\UU^{\la}_1:=\{\kkk\in \UU^{\la}:|k_1|\ge 2C_0\sqrt{n}>|k_2|\ge\cdots\ge |k_m|      \}
\end{equation*}
\begin{equation*}
\UU^{\la}_2:=\{\kkk\in \UU^{\la}:|k_1|\ge |k_2|\ge 2C_0\sqrt{n}>|k_3|\ge\cdots\ge |k_m|      \}
\end{equation*}
$$\vdots$$
\begin{equation*}
\UU^{\la}_m:=\{\kkk\in \UU^{\la}:|k_1|\ge \cdots\ge |k_m|\ge 2C_0\sqrt{n}      \}
\end{equation*}
and accordingly, 
\begin{equation*}
\si_{\ga}^{(1)}(\xxxi)=\sum_{l=1}^{m}\si_{\ga,l}^{(1)}(\xxxi)
\end{equation*}
where
\begin{equation*}
\si_{\ga,l}^{(1)}(\xxxi):={\sum_{\la: 2^{\la}\ge 2^8C_0m\sqrt{n}} \; \sum_{\GGG\in \II^{\la}}\sum_{\kkk\in \UU^{\la}_l}{b_{\GGG,\kkk}^{\la,\ga}\Psi_{G_1,k_1}^{\la}(\xi_1)\cdots \Psi_{G_m,k_m}^{\la}(\xi_m)} }.
\end{equation*}

Then it is enough to show that for each $1\le l\le m$
\begin{align}\label{app3reduction}
&\Big\Vert\sum_{\la: 2^{\la}\ge 2^8C_0m\sqrt{n}} \; \sum_{\GGG\in\II^{\la}}\sum_{\ga\in \bbz} \sum_{\kkk\in \UU_{l}^{\la}}  b_{\GGG,\kkk}^{\la,\ga} \Big(\prod_{j=1}^{m}L_{G_j,k_j}^{\la,\ga}f_j \Big)\Big\Vert_{L^{2/m}}\nonumber\\
&\lesssim \sup_{j\in\bbz}\big\Vert \si(2^j\vec{\; \cdot\;})\wh{\Phi^{(m)}}\big\Vert_{L_s^q((\bbrn)^m)}\prod_{j=1}^{m}\Vert f_j\Vert_{L^2}
\end{align}
where $L_{G,k}^{\la,\ga}$ is defined as in (\ref{lgklg}).

Observe that if $|k|\ge 2C_0\sqrt{n}$ and $|2^{\la-\ga}\xi-k|\le C_0\sqrt{n}$, then
\begin{equation*}
C_0\sqrt{n}\le |k|-C_0\sqrt{n}\le 2^{\la-\ga}|\xi|\le |k|+C_0\sqrt{n}\le 2^{\la+2}+C_0\sqrt{n}\le 2^{\la+3},
\end{equation*}
which implies
\begin{equation}\label{equiexpression}
L_{G,k}^{\la,\ga}f(x)=L_{G,k}^{\la,\ga}f^{\la,\ga}(x)
\end{equation} where $f^{\la,\ga}:=\big(\wh{f}\chi_{C_0\sqrt{n}2^{\ga-\la}\le |\cdot|\le 2^{\ga+3}}\big)^{\vee}$.
Furthermore, a direct computation with Plancherel's idendity proves
\begin{equation}\label{fll2}
\Big( \sum_{\ga\in\bbz}{\Vert f^{\la,\ga}\Vert_{L^2}^2}\Big)^{1/2}\lesssim_{C_0} (\la+3)^{1/2}\Vert f\Vert_{L^2}.
\end{equation}
Let \begin{equation*}
\mathfrak{T}_{l,\GGG}^{\la,\ga}(f_1,\dots,f_m)(x):=\sum_{\kkk\in \UU_{l}^{\la}}b_{\GGG,\kkk}^{\la,\ga} \Big(\prod_{j=1}^{l}L_{G_j,k_j}^{\la,\ga}f_j^{\la,\ga}(x)\Big)\Big( \prod_{j=l+1}^{m}L_{G_j,k_j}^{\la,\ga}f_j(x) \Big).
\end{equation*}
 Then, due to (\ref{equiexpression}), the left-hand side of (\ref{app3reduction}) is less than
\begin{equation}\label{mainkeykey}
\bigg(\sum_{\la: 2^{\la}\ge 2^8C_0m\sqrt{n}} \; \sum_{\GGG\in\II^{\la}}   \Big\Vert \sum_{\ga\in \bbz}      \mathfrak{T}_{l,\GGG}^{\la,\ga}(f_1,\dots,f_m)       \Big\Vert_{L^{2/m}}^{2/m} \bigg)^{m/2}.
\end{equation}
We claim that for $1\le l\le m$ there exists a constant $C>0$ such that
\begin{align}\label{mmaingoal}
\Big\Vert \sum_{\ga\in \bbz}      \mathfrak{T}_{l,\GGG}^{\la,\ga}(f_1,\dots,f_m)  \Big\Vert_{L^{2/m}}&\le C 2^{-\la(s-\max{( \frac{(m-1)n}{2}, \frac{mn}{q} )})}(\la+3)^{m} \\
 &~ \times \sup_{j\in\bbz}\big\Vert \si(2^j\vec{\; \cdot\;})\wh{\Phi^{(m)}}\big\Vert_{L_s^q((\bbrn)^m)}\prod_{j=1}^{m}\Vert f_j\Vert_{L^2},\nonumber
\end{align}
which clearly implies that (\ref{mainkeykey}) is majorized by the right-hand side of (\ref{app3reduction}) as 
\begin{equation*}
 \Big(  \sum_{\la: 2^{\la}\ge 2^8C_0m\sqrt{n}}   2^{-\frac{2\la}{m}(s-\max{(\frac{(m-1)n}{2},\frac{mn}{q} )})}  (\la+3)^2   \Big)^{m/2}<\infty,
\end{equation*}
which is due to the assumption $s>\max{(    \frac{(m-1)n}{2}, \frac{mn}{q}  )}$.

Therefore, let us prove (\ref{mmaingoal}).

\subsection{The case $l=1$}
We utilize the Littlewood-Paley theory for Hardy spaces as in Section \ref{pfapp2}.
 There exists a unique polynomial $Q^{\la,\GGG}(x)$ such that
\begin{align}\label{l=1}
\Big\Vert \sum_{\ga\in\bbz}  &\mathfrak{T}_{1,\GGG}^{\la,\ga}(f_1,\dots,f_m)  -Q^{\la,\GGG}  \Big\Vert_{L^{2/m}}\notag\\
\lesssim &\Big\Vert  \Big\{ \Phi_j^{(m)}\ast \Big( \sum_{\ga\in\bbz}  \mathfrak{T}_{1,\GGG}^{\la,\ga}(f_1,\dots,f_m) \Big)  \Big\}_{j\in\bbz}\Big\Vert_{L^{2/m}(\ell^2)}.
\end{align} 
Note that
\begin{equation*}
2^{\la-3}\le 2^{\la-2}-2C_0\sqrt{mn}\le |\kkk|- (|k_2|^2+\dots+|k_m|^2)^{1/2}\le |k_1|\le |\kkk|\le 2^{\la+2}
\end{equation*}
and this proves that 
\begin{equation*}
\textup{Supp}\big(\Psi_{G_1,k_1}^{\la}(\cdot/2^{\ga})\big)\subset \{\xi\in \bbrn: 2^{\ga-4}\le |\xi|\le 2^{\ga+3}\}.
\end{equation*}
Moreover, since $|k_j|\le 2C_0\sqrt{n}$ for $2\le j\le m$,
\begin{equation*}
\textup{Supp}\big(\Psi_{G_j,k_j}^{\la}(\cdot/2^{\ga})\big)\subset \{\xi\in \bbrn: |\xi|\le 2^{-6}m^{-1}2^{\ga}\}.
\end{equation*}
 Therefore, the Fourier transform of $\mathfrak{T}_{1,\GGG}^{\la,\ga}(f_1,\dots,f_m)$ is supported in the set $\{\xi\in \bbrn: 2^{\ga-5}\le |\xi|\le 2^{\ga+4} \}$ and the technique of (\ref{marshallest}) yields that the right-hand side of (\ref{l=1}) is dominated by a constant times
\begin{equation*}
 \Big\Vert   \Big( \sum_{\ga\in\bbz}\big| \mathfrak{T}_{1,\GGG}^{\la,\ga}(f_1,\dots,f_m)\big|^2 \Big)^{1/2} \Big\Vert_{L^{2/m}}.
\end{equation*}
The $L^{2/m}$-norm is bounded  by
\begin{equation*}
\sup_{\ga\in\bbz}\big\Vert \{ b_{\GGG,\kkk}^{\la,\ga}\}_{\kkk\in (\bbzn)^m}\big\Vert_{\ell^{\infty}}2^{\la mn/2}\Big( \sum_{\ga\in\bbz}\Vert  f_1^{\la,\ga}\Vert_{L^2}^2\Big)^{1/2}\prod_{j=2}^{m}\Vert f_j\Vert_{L^2}
\end{equation*}
thanks to Proposition \ref{keyapplication1}.
The embedding $\ell^{q}\hookrightarrow \ell^{\infty}$ and (\ref{lrestimate}) imply
\begin{equation}\label{suplinfty}
\sup_{\ga\in\bbz}\big\Vert \{b_{\GGG,\kkk}^{\la,\ga}\}_{\kkk\in (\bbzn)^m}\big\Vert_{\ell^{\infty}}\lesssim 2^{-\la(s-mn/q+mn/2)}\sup_{j\in\bbz}\big\Vert \si(2^j\vec{\; \cdot\;})\wh{\Phi^{(m)}}\big\Vert_{L_s^q((\bbrn)^m)}.
\end{equation}
This, together with (\ref{fll2}), finally proves that the left-hand side of (\ref{l=1}) is dominated by a constant multiple of 
\begin{equation*}
 2^{-\la(s-mn/q)}(\la+3)^{1/2}\sup_{j\in\bbz}\big\Vert \si(2^j\vec{\; \cdot\;})\wh{\Phi^{(m)}}\big\Vert_{L^q_s((\bbrn)^m)}\prod_{j=1}^{m}\Vert f_j\Vert_{L^2}
\end{equation*}
and accordingly,
\begin{equation*}
\sum_{\ga\in\bbz}  \mathfrak{T}_{1,\GGG}^{\la,\ga}(f_1,\dots,f_m)  -Q^{\la,\GGG}\in L^{2/m}.
\end{equation*}

Moreover, Proposition \ref{keyapplication1}, together with (\ref{suplinfty}), yields that
\begin{align*}
\Big\Vert \sum_{\ga\in\bbz}  \mathfrak{T}_{1,\GGG}^{\la,\ga}(f_1,\dots,f_m)   \Big\Vert_{L^{2/m}}&\le \Big\Vert \sum_{\ga\in\bbz}  \big| \mathfrak{T}_{1,\GGG}^{\la,\ga}(f_1,\dots,f_m)  \big| \Big\Vert_{L^{2/m}}\\
&\lesssim 2^{-\la(s-mn/q)}\Big( \sum_{\ga\in\mathbb{Z}}{\Vert f_1^{\la,\ga}\Vert_{L^2}}\Big)\prod_{j=2}^{m}\Vert f_j\Vert_{L^2}
\end{align*}
and, similarly to (\ref{l=1case}), we have
\begin{equation*}
\Vert f_1^{\la,\ga}\Vert_{L^2}=\bigg[ \int_{C_0\sqrt{n} 2^{\ga-\la}\le |\xi|\le 2^{\ga+3}}{|\wh{f_1}(\xi)|^2}d\xi\bigg]^{\f 12}\lesssim_{N}\begin{cases}
2^{(\ga+3)n/2},~& \ga<0\\
2^{-(\ga-\la)(N-n/2)}~& \ga\ge 0
\end{cases} 
\end{equation*} for sufficiently large $N>n/2$.
Using the argument that led to (\ref{casel=1}), we see that
\begin{equation*}
 \sum_{\ga\in\bbz}  \mathfrak{T}_{1,\GGG}^{\la,\ga}(f_1,\dots,f_m)\in L^{2/m}
\end{equation*} 
and thus $Q^{\la,\GGG}=0$. 
 Then the inequality (\ref{mmaingoal}) for $l=1$ follows.

\subsection{ The case $2\le l\le m$} 
If $0<q<\frac{2l}{l-1}$, we simply apply Proposition \ref{keyapplication1} to have
\begin{align*}
&\Big\Vert \sum_{\ga\in\bbz}{\mathfrak{T}_{l,\GGG}^{\la,\ga}(f_1,\dots,f_m)}\Big\Vert_{L^{2/m}}\\
&\lesssim \sup_{\ga\in\bbz}
\big\Vert \{b_{\GGG,\kkk}^{\la,\ga} \}_{\kkk\in (\bbzn)^m}\big\Vert_{\ell^q} 2^{\la mn/2} \Big[\prod_{j=1}^{l}\Big(\sum_{\ga\in\bbz}{\Vert f_j^{\la,\ga}\Vert_{L^2}^2} \Big)^{1/2} \Big]\Big[ \prod_{j=l+1}^{m}\Vert f_j\Vert_{L^2}\Big]
\end{align*}
where the embedding $\ell^q\hookrightarrow \ell^{\infty}$ is applied.
Then the last expression is no more than a constant multiple of
\begin{equation*}
2^{-\la(s-mn/q)}(\la+3)^{l/2}\sup_{j\in\bbz}\big\Vert \si(2^j\vec{\; \cdot\;})\wh{\Phi^{(m)}}\big\Vert_{L_s^q((\bbrn)^m)}\prod_{j=1}^{m}\Vert f_j\Vert_{L^2}
\end{equation*} by using (\ref{lrestimate}) and (\ref{fll2}). Then (\ref{mmaingoal}) follows.

If $\frac{2l}{l-1}\le q<\infty$,  using the argument in proving (\ref{mainptest}), we obtain
 \begin{align*}
&\big|\mathfrak{T}_{l,\GGG}^{\la,\ga}\big(f_1,\dots,f_m\big)(x) \big|\\
&\lesssim 2^{\la (m-l)n/2}\Big( \prod_{j=l+1}^{m}\mathcal{M}f_{j}(x)\Big)\\
&\q\times \sum_{\kkk^{*1,\dots,l}\in \mathcal{P}_{*1,\dots,l}\UU_{l}^{\la}}\Big|\sum_{\kkk^{1,\dots,l}\in Col_{\kkk^{*1,\dots,l}}^{\UU_{l}^{\la}}}b_{\GGG,\kkk}^{\la,\ga}\Big( \prod_{j=1}^{l}L_{G_j,k_j}^{\la}\big( f_j^{\la,\ga}(\cdot/2^{\ga})\big)(2^{\ga}x)\Big)\Big|
\end{align*} 
and then H\"older's inequality and the $L^2$ boundedness of $\mathcal{M}$ yield that
\begin{align*}
&\Big\Vert \sum_{\ga\in\bbz}{\mathfrak{T}_{l,\GGG}^{\la,\ga}(f_1,\dots,f_m)}\Big\Vert_{L^{2/m}}\lesssim 2^{\la(m-l)n/2}\Big(\prod_{j=l+1}^{m}\Vert f_j\Vert_{L^2} \Big)  \\
&  \times \sum_{\kkk^{*1,\dots,l}\in \mathcal{P}_{*1,\dots,l}\UU_{l}^{\la}}  \Big\Vert \sum_{\ga\in\bbz}  \Big|\sum_{\kkk^{1,\dots,l}\in Col_{\kkk^{*1,\dots,l}}^{\UU_{l}^{\la}}} 
\!\!\!\!\!\! b_{\GGG,\kkk}^{\la,\ga}\Big( \prod_{j=1}^{l}L_{G_j,k_j}^{\la}\big( f_j^{\la,\ga}(\cdot/2^{\ga}) \big)(2^{\ga}\cdot)\Big)\Big|       \Big\Vert_{L^{2/l}}
\end{align*}
where we used the fact $|\mathcal{P}_{*1,\dots,l}\UU_{l}^{\la}|\lesssim 1$.
Due to Proposition \ref{05071}, a change of variables, and  (\ref{lrestimate}), 
the $L^{2/l}$-norm is less than
\begin{align*}
& \Big( \sum_{\ga\in\bbz}  \Big\Vert \sum_{\kkk^{1,\dots,l}\in Col_{\kkk^{*1,\dots,l}}^{\UU_{l}^{\la}}}b_{\GGG,\kkk}^{\la,\ga}\Big( \prod_{j=1}^{l}L_{G_j,k_j}^{\la}\big( f_j^{\la,\ga}(\cdot/2^{\ga}) \big)(2^{\ga}\cdot)\Big)\Big\Vert_{L^{2/l}}^{2/l}    \Big)^{l/2}\\
 &= \Big( \sum_{\ga\in\bbz} 2^{-\ga n} \Big\Vert \sum_{\kkk^{1,\dots,l}\in Col_{\kkk^{*1,\dots,l}}^{\UU_{l}^{\la}}}b_{\GGG,\kkk}^{\la,\ga}\Big( \prod_{j=1}^{l}L_{G_j,k_j}^{\la}\big( f_j^{\la,\ga}(\cdot/2^{\ga}) \big)\Big)\Big\Vert_{L^{2/l}}^{2/l}    \Big)^{l/2}\\
 &\lesssim E_{q,l,\la}2^{-\la(s-mn/q+mn/2)}2^{\la ln/2}\sup_{j\in\bbz}\big\Vert \si(2^j\vec{\; \cdot\;})\wh{\Phi^{(m)}}\big\Vert_{L^q_s((\bbrn)^m)}\\
 &\qq \qq \qq \qq \qq\qq  \times \Big(\sum_{\ga\in\bbz}2^{-\ga n}     \prod_{j=1}^{l}\big\Vert f_j^{\la,\ga}(\cdot/2^{\ga})\big\Vert_{L^2}^{2/l}\Big)^{l/2}
 \end{align*}
 where 
\begin{equation*}
E_{q,l,\la}:=\begin{cases}
\la^{l/2}, & q=\frac{2l}{l-1}\\
2^{\la n(l/2-l/q-1/2)}, & q>\frac{2l}{l-1}
\end{cases}.
\end{equation*}
Since
\begin{align*}
& \Big( \sum_{\ga\in\bbz}  2^{-\ga n}\prod_{j=1}^{l}\big\Vert f_j^{\la,\ga}(\cdot/2^{\ga})\big\Vert_{L^2}^{2/l} \Big)^{l/2} \le \prod_{j=1}^{l}\Big( \sum_{\ga \in\bbz}{2^{-\ga n}\Vert f_j^{\la,\ga}(\cdot/2^{\ga})\Vert_{L^2}^2}\Big)^{1/2}\\
  &\qq\qq\qq\qq \le\prod_{j=1}^{l}\Big( \sum_{\ga\in\bbz}\Vert f_j^{\la,\ga}\Vert_{L^2}^2\Big)^{1/2}\lesssim (\la+3)^{l/2}\prod_{j=1}^{l}\Vert f_j\Vert_{L^2},
\end{align*}
 we finally obtain that 
\begin{equation*}
\Big\Vert \sum_{\ga\in\bbz}{\mathfrak{T}_{l,\GGG}^{\la,\ga}(f_1,\dots,f_m)}\Big\Vert_{L^{2/m}}\lesssim  F_{q,l,\la}^{(s,m,n)} \sup_{j\in\bbz}\big\Vert \si(2^j\vec{\; \cdot\;})\wh{\Phi^{(m)}}\big\Vert_{L^q_s((\bbrn)^m)}          \prod_{j=1}^{m}{\Vert f_j\Vert_{L^2}}
\end{equation*}
where
\begin{equation*}
F_{q,l,\la}^{(s,m,n)}:=E_{q,l,\la} 2^{-\la(s-mn/q)}(\la+3)^{l/2}.
\end{equation*}
It is easy to check that for $2\le l\le m$ and $\frac{2l}{l-1}\le q$
\begin{equation*}
F_{q,l,\la}^{(s,m,n)}\lesssim 2^{-\la(s-\max{(\frac{(m-1)n}{2},\frac{mn}{q}    )})}(\la+3)^m
\end{equation*}
and the proof of (\ref{mmaingoal}) is complete.

\section{Concluding remarks}

In this article  we    focused on the $L^2\times\cdots\times L^2\to L^{2/m}$ boundedness 
for several fundamental $m$-linear operators.  In future work we plan to 
obtain similar initial estimates for maximal singular integrals and 
maximal multipliers.

The $L^{2}\times \cdots \times L^{2} $ estimates obtained in this paper 
provide   crucial 
 initial   bounds  that provide the cornerstone needed to  launch a complete 
 boundedness study   on general products 
 of Lebesgue spaces. Certainly our initial estimates can be 
extended to include   points  obtained by duality and 
interpolation; these are called local $L^2$ points. 
For the remaining points there are techniques available, for instance,   
interpolation between dyadic pieces of an operator between   {\it good} 
local $L^2$ points and {\it bad} points near the boundary of the region 
$1<p_1,\dots, p_m<\nf$, $1/m<p<\nf$; this technique was developed in 
\cite{Gr_He_Ho} in the bilinear case. 
We chose not to pursue this line of investigation here in order to 
 direct our focus on the idea of wavelet expansions
 and shorten the exposition.   We plan to pursue 
general  $L^{p_1}\times \cdots \times L^{p_m}\to L^p$ 
boundedness for many multilinear  operators  
 in   subsequent work. 
It should be mentioned that in 
 a recent manuscript of Heo, Lee, Hong, Yang, Lee, and Park \cite{6author} 
 the extension to the full range of indices was obtained for Theorem~\ref{application3}, when $q=2$, although
 the case of general $q$ remains unresolved.  
 
 \medskip

 \noindent{\bf Acknowledgement.} We are indebted to Professor Miyachi for pointing out  to us the relation between  \cite{Katoetal1} and Proposition~\ref{maintheorem}.

%\bibliographystyle{apalike}
%\bibliography{HeBib}

%\printindex

\end{document}